\documentclass{amsart}

\usepackage[english]{babel}

\usepackage{amsmath}

\usepackage[mathx,matha]{mathabx}

\usepackage{amstext}

\usepackage{amssymb}

\usepackage{array}

\usepackage{amsthm}

\usepackage{dsfont}

\usepackage{dirtytalk}

\usepackage[latin1]{inputenc}

\usepackage[T1]{fontenc}

\usepackage{mathtools}

\usepackage{appendix}

\usepackage{xcolor}

\usepackage{ulem}

\usepackage[arrow, matrix, curve]{xy}

\usepackage{verbatim}

\usepackage{graphicx}

\usepackage{color}

\usepackage{transparent}

\usepackage{subfigure}

\usepackage{hyperref}

\usepackage{cases}

\usepackage[update,prepend]{epstopdf}

\usepackage{bbold}

\newtheorem{theorem}{Theorem}[section]

\newtheorem*{theorem*}{Theorem}

\theoremstyle{plain}

\newtheorem*{conjecture*}{Conjecture}

\newtheorem{assumption}[theorem]{Assumption}

\newtheorem{corollary}[theorem]{Corollary}

\newtheorem{proposition}[theorem]{Proposition}

\newtheorem{lemma}[theorem]{Lemma}

\theoremstyle{remark}

\newtheorem{condition}{Condition}

\theoremstyle{definition}

\newtheorem{definition}[theorem]{Definition}

\newtheorem{remark}[theorem]{Remark}

\def\eps{\varepsilon}

\def\bi{\begin{itemize}}

\def\ei{\end{itemize}}

\newcommand{\inner}[2]{\ensuremath{\left\langle #1 , #2 \right\rangle}}

\newcommand{\R}{\mathbb{R}}

\newcommand{\Z}{\mathbb{Z}}

\renewcommand{\dim}{{\rm dim}\,}

\newcommand{\conv}{{\rm conv}}

\renewcommand{\phi}{\varphi}

\DeclareMathOperator{\spann}{span}

\DeclareOldFontCommand{\it}{\normalfont\itshape}{\mathit}

\newcommand{\bspm}{\left(\begin{smallmatrix}}\newcommand{\espm}{\end{smallmatrix}\right)}

\newcommand{\bpm}{\begin{pmatrix}}\newcommand{\epm}{\end{pmatrix}}

\def\bs{\begin{satz}}\def\es{\end{satz}}

\def\blem{\begin{lemma}}\def\elem{\end{lemma}}

\def\bthm{\begin{theorem}}\def\ethm{\end{theorem}}

\def\bcor{\begin{corollary}}\def\ecor{\end{corollary}}

\def\beq{\begin{equation}}\def\eeq{\end{equation}}

\def\beqq{\begin{equation*}}\def\eeqq{\end{equation*}}

\def\bal{\begin{align}}\def\eal{\end{align}}

\def\bpf{\begin{proof}}\def\epf{\end{proof}}

\def\bex{\begin{example}}\def\eex{\end{example}}

\def\brem{\begin{remark}}\def\erem{\end{remark}}

\def\bass{\begin{assumption}}\def\eass{\end{assumption}}

\def\bprop{\begin{proposition}}\def\eprop{\end{proposition}}

\def\bdefi{\begin{definition}}\def\edefi{\end{definition}}

\def\bcond{\begin{condition}}\def\econd{\end{condition}}

\def\bconj{\begin{conjecture*}}\def\econj{\end{conjecture*}}

\DeclareFontEncoding{FMS}{}{}

\DeclareFontSubstitution{FMS}{futm}{m}{n}

\DeclareFontEncoding{FMX}{}{}

\DeclareFontSubstitution{FMX}{futm}{m}{n}

\DeclareSymbolFont{fouriersymbols}{FMS}{futm}{m}{n}

\DeclareSymbolFont{fourierlargesymbols}{FMX}{futm}{m}{n}

\DeclareMathDelimiter{\VERT}{\mathord}{fouriersymbols}{152}{fourierlargesymbols}{147}

\def\bi{\begin{itemize}}

\def\ei{\end{itemize}}

\def\ben{\begin{enumerate}}

\def\een{\end{enumerate}}

\usepackage{enumitem}

\usepackage{tcolorbox}

\newtcolorbox{implementation}[2][]{colframe=blue!75!black,colbacktitle=green!10!white,colback=green!10!white,coltitle=green!75!black,title={#2},fonttitle=\bfseries,#1}

\begin{document}

\title[Regularity result for shortest billiard trajectories]{A regularity result for shortest generalized billiard trajectories in convex bodies in $\R^n$}

\author{Daniel Rudolf and Stefan Krupp}

%\author{Stefan Krupp$^\blacklozenge$ and Daniel Rudolf$^\clubsuit$}

\date{\today}

\maketitle
%\tableofcontents

\begin{abstract}
We study length-minimizing closed generalized Euclidean billiard trajectories in convex bodies in $\R^n$ and investigate their relation to the inclusion minimal affine sections that contain these trajectories. We show that when passing to these sections, the length-minimizing closed billiard trajectories are still billiard trajectories, but their length-minimality as well as their regularity can be destroyed. In light of this, we prove what weaker regularity is actually preserved under passing to these sections. Based on the results, we develop an algorithm in order to calculate length-minimizing closed regular billiard trajectories in convex polytopes in $\R^n$.
\end{abstract}

\section{Introduction}\label{Sec:Intro}

In this paper, we analyze closed billiard trajectories in convex bodies in $\R^n$. The billiard trajectories are meant to be \textit{Euclidean} what for bouncing points locally means: the angle of reflection equals the angle of incidence (here, we assume that normal vector as well as incident and reflected rays lie in the same two-dimensional affine flat). This local reflection rule can be seen as consequence of the global least action principle. For a reflection in a hyperplane this principle means that a billiard trajectory segment $(p_{j-1},p_j,p_{j+1})$ minimizes the Euclidean length in the space of all paths connecting $p_{j-1}$ and $p_{j+1}$ via a reflection at this hyperplane.

From the geometric optics point of view, Euclidean billiards describe the wave propagation in a medium which is not only homogeneous and isotropic but also contains perfectly reflecting mirrors.

There is generally much interest in the study of billiards: Problems in various mathematical fields can be related to problems in mathematical billiards; see for example \cite{Gutkin2012}, \cite{Katok2005}, and \cite{Tabachnikov2005} for comprehensive surveys. Euclidean billiard trajectories in the plane have been investigated intensively. Nonetheless, so far, not much is known about Euclidean billiard trajectories on higher-dimensional \say{tables}.

Let us precisely define closed Euclidean billiard trajectories on convex bodies in $\R^n$ based on the above mentioned least action principle.

\begin{figure}[h!]
\centering
\def\svgwidth{200pt}
\begingroup%
  \makeatletter%
  \providecommand\color[2][]{%
    \errmessage{(Inkscape) Color is used for the text in Inkscape, but the package 'color.sty' is not loaded}%
    \renewcommand\color[2][]{}%
  }%
  \providecommand\transparent[1]{%
    \errmessage{(Inkscape) Transparency is used (non-zero) for the text in Inkscape, but the package 'transparent.sty' is not loaded}%
    \renewcommand\transparent[1]{}%
  }%
  \providecommand\rotatebox[2]{#2}%
  \newcommand*\fsize{\dimexpr\f@size pt\relax}%
  \newcommand*\lineheight[1]{\fontsize{\fsize}{#1\fsize}\selectfont}%
  \ifx\svgwidth\undefined%
    \setlength{\unitlength}{285.81452862bp}%
    \ifx\svgscale\undefined%
      \relax%
    \else%
      \setlength{\unitlength}{\unitlength * \real{\svgscale}}%
    \fi%
  \else%
    \setlength{\unitlength}{\svgwidth}%
  \fi%
  \global\let\svgwidth\undefined%
  \global\let\svgscale\undefined%
  \makeatother%
  \begin{picture}(1,0.69875334)%
    \lineheight{1}%
    \setlength\tabcolsep{0pt}%
    \put(0,0){\includegraphics[width=\unitlength,page=1]{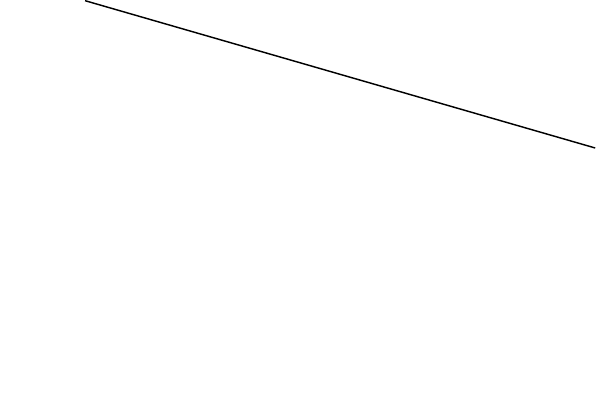}}%
    \put(0.86349328,0.42473384){\color[rgb]{0,0,0}\makebox(0,0)[lt]{\lineheight{1.25}\smash{\begin{tabular}[t]{l}$H_j$\end{tabular}}}}%
    \put(0.71786167,0.0082515){\color[rgb]{0,0,0}\makebox(0,0)[lt]{\lineheight{1.25}\smash{\begin{tabular}[t]{l}$p_{j+1}$\end{tabular}}}}%
    \put(-0.00290425,0.09233721){\color[rgb]{0,0,0}\makebox(0,0)[lt]{\lineheight{1.25}\smash{\begin{tabular}[t]{l}$p_{j-1}$\end{tabular}}}}%
    \put(0.56097889,0.60368355){\color[rgb]{0,0,0}\makebox(0,0)[lt]{\lineheight{1.25}\smash{\begin{tabular}[t]{l}$p_j$\end{tabular}}}}%
    \put(0.0728596,0.42675179){\color[rgb]{0,0,0}\makebox(0,0)[lt]{\lineheight{1.25}\smash{\begin{tabular}[t]{l}$T$\end{tabular}}}}%
    \put(0,0){\includegraphics[width=\unitlength,page=2]{Stossregelmink.pdf}}%
    \put(0.38372694,0.64717047){\color[rgb]{0,0,0}\makebox(0,0)[lt]{\lineheight{1.25}\smash{\begin{tabular}[t]{l}$\widebar{p}_j$\end{tabular}}}}%
    \put(0,0){\includegraphics[width=\unitlength,page=3]{Stossregelmink.pdf}}%
  \end{picture}%
\endgroup%
\caption{The billiard reflection rule: $p_j$ minimizes \eqref{eq:minimization} over all $\widebar{p}_j\in H_j$, where $H_j$ is a $T$-supporting hyperplane through $p_j$.}
\label{img:Stossregel}
\end{figure}

\bdefi\label{Def:closedbilliardtrajectory}
Let $T\subset\R^n$ be a convex body, i.e., a compact convex set with non-empty interior, which from now on we call the \textit{billiard table}. We say that a closed polygonal curve\footnote{For the sake of simplicity, whenever we talk of the vertices $p_1,...,p_m$ of a closed polygonal curve, we assume that they satisfy $p_j\neq p_{j+1}$ and $p_j$ is not contained in the line segment connecting $p_{j-1}$ and $p_{j+1}$ for all $j\in\{1,...,m\}$. Furthermore, whenever we settle indices $1,...,m$, then the indices in $\Z$ will be considered as indices modulo $m$.} with vertices $p_1,...,p_m$, $m\geq 2$, on the boundary of $T$ (denoted by $\partial T$) is a \textit{closed billiard trajectory on $T$} if for every $j\in \{1,...,m\}$ there is a $T$-supporting hyperplane $H_j$ through $p_j$ such that $p_j$ minimizes
\beq ||\widebar{p}_j-p_{j-1}||+||p_{j+1}-\widebar{p}_j||\label{eq:minimization}\eeq
over all $\widebar{p}_j\in H_j$. We encode this closed billiard trajectory by $( p_1,...,p_m )$ and call its vertices \textit{bouncing points}. Its \textit{length}\footnote{In what follows, this is the way how we measure the length of closed polygonal curves in general.} is given by
\beq \ell(( p_1,...,p_m))=\sum_{j=1}^m ||p_{j+1}-p_j||.\label{eq:length}\eeq
\edefi

We call a boundary point $p_j\in \partial T$ \textit{smooth} if there is a unique $T$-supporting hyperplane through $p_j$. We say that $\partial T$ is \textit{smooth} if every boundary point is smooth. We call a closed billiard trajectory \textit{regular} if all its bouncing points are smooth boundary points of the billiard table boundary.

We remark that the notion of billiard trajectories is usually used for classical trajectories, i.e., for trajectories with bouncing points in smooth boundary points (billiard table cushions) while they terminate in non-smooth boundary points (billiard table pockets). Definition \ref{Def:closedbilliardtrajectory} generalizes this classical billiard reflection rule to non-smooth boundaries. To the authors' knowledge, the papers \cite{Bezdek1989} ('89), \cite{Ghomi2004} ('04), and \cite{Bezdek2011} ('09) were among the first suggesting a detailed study of these \textit{generalized} billiard trajectories.

Definition \ref{Def:closedbilliardtrajectory} implies the local billiard reflection rule: the hyperplane associated to the corresponding normal vector within the local formulation is the one appearing within the minimization of \eqref{eq:minimization}. Since a $T$-supporting hyperplane through non-smooth boundary points of the billiard table $T$ is not unique, from a constructive point of view, the billiard reflection rule may produce different bouncing points following two already known consecutive ones.

The main interest we pursue in the first part of this paper is the investigation of length-minimizing closed billiard trajectories on a convex body $T\subset\R^n$ in relation to the inclusion minimal affine sections that contain these trajectories. Its relevance is given against the background of the development of an algorithm for computing length-minimizing closed billiard trajectories which we will discuss in the second part of this paper.

In a nutshell, our result can be summarized as follows: We show that when passing to the inclusion minimal affine section that contains the length-minimizing billiard trajectory, this length-minimizing billiard trajectory is still a billiard trajectory, but almost every other property (e.g., length-minimality, regularity) can be destroyed (see examples (B) and (C) in Section \ref{Sec:Examples}). In light of this observation, we prove (see Theorem \ref{Thm:RegularityResult1}) what weaker regularity is actually preserved under passing to the inclusion minimal affine section. Briefly, this weaker regularity can be described as follows: for a length-minimizing billiard trajectory the billiard reflections are uniquely determined if restricted to the inclusion minimal affine section that contains this trajectory; that is, the billiard reflection rule is unambigious provided that the reflected ray lies in the inclusion minimal affine section.

In order to state this weaker regularity result, let
\beqq N_T(p):=\{n\in\R^n : \langle n,y-p\rangle \leq 0 \text{ for all }y\in T\}\eeqq
be the \textit{outer normal cone} of $T$ at the point $p\in \partial T$. Then, our result reads as follows:

\bthm\label{Thm:RegularityResult1}
Let $T\subset\R^n$ be a billiard table and $p=(p_1,...,p_m)$ a length-minimizing closed billiard trajectory on $T$. Further, let $U\subseteq\R^n$ be the convex cone spanned by the normal vectors related to the billiard reflection rule and let $V\subseteq\R^n$ be the affine subspace such that $T\cap V$ is the inclusion minimal affine section of $T$ that contains $p$. Then, it follows that
\beqq U=V_0,\quad \dim U = \dim V = m-1\eeqq
and
\beq \dim \left(N_T(p_j)\cap V_0\right)=1\label{eq:RegularityResult1}\eeq
for all $j\in\{1,...,m\}$, where $V_0$ is the linear subspace that is parallel to the affine subspace $V$ satisfying $\dim V_0 = \dim V$.
\ethm

From
\beqq \dim U = \dim V\leq n,\eeqq
it follows that
\beqq m\leq n+1.\eeqq
Furthermore,
\beqq \dim V=m-1\eeqq
implies that $p$ is \textit{maximally spanning} by what we mean
\beqq \dim (\conv\{p_1,...,p_m\})=m-1.\eeqq

In fact, \eqref{eq:RegularityResult1} is a regularity result: If $m=n+1$, meaning that $V=V_0=\R^n$, then \eqref{eq:RegularityResult1} becomes
\beqq \dim(N_T(p_j))=1\quad \forall j\in\{1,...,m\},\eeqq
i.e., $p$ is \textit{regular}. For $n=2$ this means that every length-minimizing closed billiard trajectory on $T$ has either two or three bouncing points, while in the latter case all of them are smooth boundary points of $T$.

Some special cases of Theorem \ref{Thm:RegularityResult1} were already known: For length-minimizing closed billiard trajectories it has been proven in \cite{Bezdek2011} that $m$ is bounded from above by $n+1$. The two-dimensional case has been shown in \cite{AlkoumiSchlenk2014}, while in \cite{AkopBal2015}, it has been proven that every length-minimizing closed billiard trajectory with $n+1$ bouncing points is regular.

We remark that in our upcoming paper \cite{KruppRudolf2022} we extend Theorem \ref{Thm:RegularityResult1} from the Euclidean setting to the more general Finsler/Minkowski setting.

The relevance of Theorem \ref{Thm:RegularityResult1} lies especially in its sharpness (which is ensured by example (C) in Section \ref{Sec:Examples}): It refutes the presumption, which at first appears to be intuitively correct, that every length-minimizing closed billiard trajectory with more than two bouncing points is regular within the inclusion minimal affine section of $T$ that contains this trajectory, i.e., that it is \textit{relatively regular} with respect to this inclusion minimal affine section. The latter could have been expected as higher dimensional generalization of the two-dimensional regularity result in \cite{AlkoumiSchlenk2014}.

The algorithm which we develop in the second part of this paper can be seen as one of the first steps in a direction towards tackling certain originally geometrical/dynamical problems arising in the field of symplectic geometry/dynamics computationally. This results from the following relationships: As was already pointed out in \cite{AlkoumiSchlenk2014}, the search for length-minimizing closed generalized billiard trajectories is related to the Ekeland-Hofer-Zehnder capacity of the $4$-dimensional Lagrangian product
\beqq T\times B^2_1(0)\eeqq
in the standard symplectic space $\R^4$, where we denote by $B^2_1(0)\subset\R^2$ the Euclidean unit ball. In relation to the volume, the aforementioned in turn is related to so far not solved isoperimetric-like/systolic billiard inequalities (see \cite{Balitskiy2018}) as well as to certain configurations in light of Viterbo's conjecture from symplectic geometry (see \cite{Viterbo2000}), which, due to \cite{ArtKarOst2013}, is related to the famous Mahler conjecture from convex geometry (see \cite{Mahler1939}). Another interesting relationship is with Wetzel's problem from geometry (see \cite{Wetzel1973}), on which we will elaborate in an upcoming paper in more detail (see \cite{Rudolf2022}).

Let us briefly present the structure of this paper: In Section \ref{Sec:PropertiesBilliard}, we discuss properties of closed generalized billiard trajectories in convex bodies in $\R^n$. In Section \ref{Sec:Proof1}, we prove Theorem \ref{Thm:RegularityResult1}. In Section \ref{Sec:Examples}, we present relevant examples which, i.a., provide the sharpness of Theorem \ref{Thm:RegularityResult1}; and, finally, in Section \ref{Sec:Construction}, we discuss the algorithm in order to calculate length-minimizing closed regular billiard trajectories in convex polytopes in $\R^n$.

\section{Properties of closed billiard trajectories}\label{Sec:PropertiesBilliard}

For a convex body $T\subset\R^n$, we define
\beqq F(T):=\left\{F\subset\R^n: F \text{ cannot be translated into }\mathring{T}\right\},\eeqq
where by $\mathring{T}$ we denote the interior of $T$ in $\R^n$ and by $\mathcal{P}(\R^n)$ the power set of $\R^n$. Based on that, we recall the following lemma from \cite{Bezdek2011}:

\blem[Lemmata 2.1, 2.2, and 2.3 in \cite{Bezdek2011}]\label{Lem:Bezdek}
\begin{itemize}
\item[(i)] Let $T\subset\R^n$ be a convex set and $F$ a finite set of at least $n+1$ points in $\R^n$. Then, there is a translate of $T$ that covers $F$ if and only if every selection of $n+1$ points of $F$ can be covered by a translate of $T$.
\item[(ii)] Let $T\subset\R^n$ be a convex body and $F$ a set of points $p_1,...,p_m\in\R^n$. Then, $F\in F(T)$ is equivalent to: there is a translate $T'$ of $T$, a selection $\{i_1,...,i_{\widetilde{m}}\}\subseteq \{1,...,m\}$, and there are closed half-spaces $H_{i_1}^+,...,H_{i_{\widetilde{m}}}^+$ of $\R^n$ such that
\begin{itemize}
\item[\tiny$\bullet$] $p_{i_j}\in \partial H_{i_j}^+$ for all $j\in\{1,...,\widetilde{m}\}$,
\item[\tiny$\bullet$] $T'\subseteq H_{i_j}^+$ for all $j\in\{1,...,\widetilde{m}\}$,
\item[\tiny$\bullet$] $\bigcap_{j=1}^{\widetilde{m}} H_{i_j}^+$ is nearly bounded, i.e. it lies between two parallel hyperplanes.
\end{itemize}
\item[(iii)] Let $T\subset\R^n$ be a billiard table and $p=(p_1,...,p_m)$ a closed billiard trajectory on $T$. Then, we have $p\in F(T)$.
\end{itemize}
\elem

We note that Lemma \ref{Lem:Bezdek}(i) can be equivalently (by contra-position applied to $\mathring{T}$) expressed by: a set $F$ of at least $n+1$ points in $\R^n$ is in $F(T)$ if and only if there is a selection of $n+1$ points out of $F$ that is in $F(T)$.

Furthermore, we recall the following characterization of length-minimizing closed billiard trajectories:

\bthm [Theorem 1.1 and Lemma 2.4 in \cite{Bezdek2011}]\label{Thm:Bezdek}
Let $T\subset\R^n$ be a billiard table. Then, every length-minimizing closed billiard trajectory on $T$ has at most $n+1$ bouncing points. Moreover, every length-minimizing closed billiard trajectory on $T$ is a length-minimizing closed polygonal curve in $F(T)$ and, conversely, every length-minimizing closed polygonal curve in $F(T)$ can be translated to a length-minimizing closed billiard trajectory on $T$.
\ethm

We note that Theorem \ref{Thm:Bezdek} is an existence result for closed billiard trajectories on arbitrary billiard tables: Let $T\subset\R^n$ be any billiard table. Then, there is a closed billiard trajectory on $T$ (with at most $n+1$ bouncing points). This can be easily concluded by a compactness argument applied to the set of closed polygonal curves in $F(T)$ combined with the $d_H$-continuity, i.e., continuity with respect to the Hausdorff distance $d_H$, of the length-functional.

We continue by stating the following property of closed billiard trajectories:

\bprop\label{Prop:SectionInvariance}
Let $T\subset\R^n$ be a billiard table, $p=(p_1,...,p_m)$ a closed billiard trajectory on $T$, and $V\subseteq\R^n$ an affine subspace such that $T\cap V$ is an affine section of $T$ that contains $p$. Then, $p$ is a closed billiard trajectory in $T\cap V$.
\eprop

\bpf
Since $p=(p_1,...,p_m)$ is a closed billiard trajectory on $T$, there are $T$-supporting hyperplanes $H_1,...,H_m$ in $\R^n$ through $p_1,...,p_m$ such that $p_j$ minimizes
\beq ||\widebar{p}_j-p_{j-1}||+||p_{j+1}-\widebar{p}_j||\label{eq:minimizationprop}\eeq
over all $\widebar{p}_j\in H_j$ for all $j\in\{1,...,m\}$. Since $T\cap V$ contains $p$, it follows that $p_j$ minimizes \eqref{eq:minimizationprop} over all $\widebar{p}_j\in H_j\cap V$ for all $j\in\{1,...,m\}$. This implies that $p$ is a billiard trajectory in $T\cap V$.
\epf

Clearly, the converse is not true: We can imagine an affine subspace $V\subseteq\R^n$ such that the section $T\cap V$ of $T$ can be translated into $\mathring{T}$. Then, every closed billiard trajectory $p$ in $T\cap V$ can be translated into $\mathring{T}$. But by Lemma \ref{Lem:Bezdek}(iii), then, $p$ cannot be a closed billiard trajectory on $T$.

In Section \ref{Sec:Examples}, we will see that, generally, the length-minimality of a closed billiard trajectory on $T$ is not invariant under going to (inclusion minimal) affine sections of $T$ that contain the closed billiard trajectory.

For what follows, it will be useful to reformulate the billiard reflection rule in the sense of the following Proposition \ref{Prop:System}. For that, we denote the $(n-1)$-dimensional unit sphere of $\R^n$ by $S^{n-1}$.

\bprop\label{Prop:System}
Let $T\subset\R^n$ be a billiard table. A closed polygonal curve with vertices $p_1,...,p_m$ on $\partial T$ is a closed billiard trajectory on $T$ if and only if there are vectors $n_1,...,n_m\in S^{n-1}$ such that
\beq \begin{cases}p_{j+1}-p_j = \lambda_j n_j, \; \lambda_j>0, \\ n_{j+1}-n_j \in -N_T(p_{j+1}) \end{cases}\label{eq:System}\eeq
is satisfied for all $j\in\{1,...,m\}$.
\eprop

\begin{figure}[h!]
\centering
\def\svgwidth{350pt}
\begingroup%
  \makeatletter%
  \providecommand\color[2][]{%
    \errmessage{(Inkscape) Color is used for the text in Inkscape, but the package 'color.sty' is not loaded}%
    \renewcommand\color[2][]{}%
  }%
  \providecommand\transparent[1]{%
    \errmessage{(Inkscape) Transparency is used (non-zero) for the text in Inkscape, but the package 'transparent.sty' is not loaded}%
    \renewcommand\transparent[1]{}%
  }%
  \providecommand\rotatebox[2]{#2}%
  \newcommand*\fsize{\dimexpr\f@size pt\relax}%
  \newcommand*\lineheight[1]{\fontsize{\fsize}{#1\fsize}\selectfont}%
  \ifx\svgwidth\undefined%
    \setlength{\unitlength}{432.27079447bp}%
    \ifx\svgscale\undefined%
      \relax%
    \else%
      \setlength{\unitlength}{\unitlength * \real{\svgscale}}%
    \fi%
  \else%
    \setlength{\unitlength}{\svgwidth}%
  \fi%
  \global\let\svgwidth\undefined%
  \global\let\svgscale\undefined%
  \makeatother%
  \begin{picture}(1,0.46554084)%
    \lineheight{1}%
    \setlength\tabcolsep{0pt}%
    \put(0,0){\includegraphics[width=\unitlength,page=1]{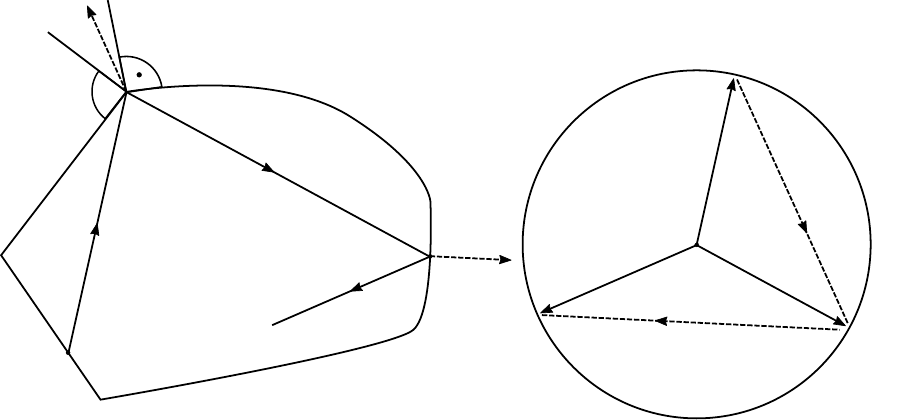}}%
    \put(0.59883995,0.35488726){\color[rgb]{0,0,0}\makebox(0,0)[lt]{\lineheight{1.25}\smash{\begin{tabular}[t]{l}$S^{n-1}$\end{tabular}}}}%
    \put(0.81036564,0.39663059){\color[rgb]{0,0,0}\makebox(0,0)[lt]{\lineheight{1.25}\smash{\begin{tabular}[t]{l}$n_{j-1}$\end{tabular}}}}%
    \put(0.9551955,0.09597517){\color[rgb]{0,0,0}\makebox(0,0)[lt]{\lineheight{1.25}\smash{\begin{tabular}[t]{l}$n_{j}$\end{tabular}}}}%
    \put(0.52974418,0.10507558){\color[rgb]{0,0,0}\makebox(0,0)[lt]{\lineheight{1.25}\smash{\begin{tabular}[t]{l}$n_{j+1}$\end{tabular}}}}%
    \put(0.01993222,0.05996557){\color[rgb]{0,0,0}\makebox(0,0)[lt]{\lineheight{1.25}\smash{\begin{tabular}[t]{l}$p_{j-1}$\end{tabular}}}}%
    \put(0.4884735,0.15723257){\color[rgb]{0,0,0}\makebox(0,0)[lt]{\lineheight{1.25}\smash{\begin{tabular}[t]{l}$p_{j+1}$\end{tabular}}}}%
    \put(0.14580639,0.33265616){\color[rgb]{0,0,0}\makebox(0,0)[lt]{\lineheight{1.25}\smash{\begin{tabular}[t]{l}$p_j$\end{tabular}}}}%
    \put(0.00538708,0.37222576){\color[rgb]{0,0,0}\makebox(0,0)[lt]{\lineheight{1.25}\smash{\begin{tabular}[t]{l}$N_T(p_j)$\end{tabular}}}}%
    \put(0.33013671,0.32261762){\color[rgb]{0,0,0}\makebox(0,0)[lt]{\lineheight{1.25}\smash{\begin{tabular}[t]{l}$T$\end{tabular}}}}%
    \put(0,0){\includegraphics[width=\unitlength,page=2]{System.pdf}}%
  \end{picture}%
\endgroup%
\caption{The visualization of \eqref{eq:System}.}
\label{img:System}
\end{figure}

\bpf
Let $p=(p_1,...,p_m)$ be a closed polygonal curve with vertices on $\partial T$. Let us assume there are vectors $n_1,...,n_m\in S^{n-1}$ which together with $p$ satisfy \eqref{eq:System}. Then, for all $j\in\{1,...,m\}$ there is a unit vector $n_T(p_{j+1})\in N_T(p_{j+1})$ such that
\beqq \begin{cases}p_{j+1}-p_j=\lambda_j n_j,\; \lambda_j >0,\\ n_{j+1}-n_j=-\mu_{j+1}n_T(p_{j+1}),\; \mu_{j+1}>0,\end{cases}\eeqq
holds. We define $H_1,...,H_m$ to be the $T$-supporting hyperplanes in $\R^n$ through $p_1,...,p_m$ which are normal to $n_T(p_1),...,n_T(p_m)$. Then, the following holds for all $j\in\{1,...,m\}$:
\allowdisplaybreaks{\begin{align*}
\nabla_{\widebar{p}_j=p_j}(||\widebar{p}_j-p_{j-1}||+||p_{j+1}-\widebar{p}_j||) = &  \frac{p_j-p_{j-1}}{||p_j-p_{j-1}||}-\frac{p_{j+1}-p_{j}}{||p_{j+1}-p_{j}||}\\
=&\frac{\lambda_{j-1}n_{j-1}}{||\lambda_{j-1}n_{j-1}||} - \frac{\lambda_j n_j}{||\lambda_j n_j||}\\
=&n_{j-1}-n_j\\
=&\mu_j n_T(p_j).
\end{align*}}%
Therefore, for all $j\in\{1,...,m\}$, $p_j$ extremizes
\beq ||\widebar{p}_j-p_{j-1}||+||p_{j+1}-\widebar{p}_j||\label{eq:system2}\eeq
over all $\widebar{p}_j\in H_j$ near $p_j$. Since for all $j\in\{1,...,m\}$, \eqref{eq:system2} is a convex function with respect to $\widebar{p}_j$, it follows for all $j\in\{1,...,m\}$ that $p_j$ is a global minimizer of \eqref{eq:system2} over all $\widebar{p}_j\in H_j$. Therefore, the billiard reflection rule is satisfied in $p_j$ for all $j\in\{1,...,m\}$. Eventually, $p$ is a closed billiard trajectory on $T$.

Conversely, let us assume $p=(p_1,...,p_m)$ is a closed billiard trajectory on $T$. Then, there are $T$-supporting hyperplanes $H_1,...,H_m$ in $\R^n$ through $p_1,...,p_m$ such that for all $j\in\{1,...,m\}$, $p_j$ minimizes \eqref{eq:system2} over all $\widebar{p}_j\in H_j$. By Lagrange's multiplier theorem, this means
\beqq \frac{p_j-p_{j-1}}{||p_j-p_{j-1}||}-\frac{p_{j+1}-p_{j}}{||p_{j+1}-p_{j}||} = \mu_j n_T(p_j),\; \mu_j\in\R,\eeqq
where $n_T(p_j)$ is the outer unit vector normal to $H_j$, for all $j\in\{1,...,m\}$. Taking inner product with $n_T(p_j)$ gives $\mu_j>0$ for all $j\in\{1,...,m\}$. If we define
\beq n_j:=(p_{j+1}-p_j)/\lambda_j,\quad \lambda_j:=||p_{j+1}-p_j||,\quad  j\in\{1,...,m\},\label{eq:system3}\eeq
and consider $n_T(p_j)\in N_T(p_j)$ for all $j\in \{1,...,m\}$, then \eqref{eq:System} is satisfied for $p$ together with the unit vectors $n_1,...,n_m$ defined in \eqref{eq:system3}.
\epf

The proof of Proposition \ref{Prop:System} shows even more: A closed polygonal curve with vertices $p_1,...,p_m$ on $\partial T$ is a closed billiard trajectory on $T$ with $H_1,...,H_m$ as $T$-supporting hyperplanes that are associated to the billiard reflection rule if and only if there are vectors $n_1,...,n_m\in S^{n-1}$ such that
\beqq \begin{cases}p_{j+1}-p_j = \lambda_j n_j, \; \lambda_j>0, \\ n_{j+1}-n_j = -\mu_{j+1} n_{H_{j+1}},\; \mu_{j+1}>0, \end{cases}\eeqq
is satisfied for all $j\in\{1,...,m\}$, where we denoted the outer unit normal vectors at $H_1,...,H_m$ by $n_{H_1},...,n_{H_m}$.

The following rather obvious proposition is needed within the proof of Theorem \ref{Thm:RegularityResult1}. It follows immediately from within the proof of Proposition \ref{Prop:System}.

\bprop\label{Prop:onlyonehyperplane}
Let $T\subset\R^n$ be a billiard table and $p=(p_1,...,p_m)$ a closed billiard trajectory on $T$. Then, for every $j\in\{1,...,m\}$ there is only one $T$-supporting hyperplane through $p_j$ for which the billiard reflection rule in $p_j$ is satisfied.
\eprop

\bpf
This claim follows from the fact that the outer unit vector $n_T(p_j)$ that is normal to $H_j$ is uniquely determined by the condition
\beqq \frac{p_j-p_{j-1}}{||p_j-p_{j-1}||}-\frac{p_{j+1}-p_{j}}{||p_{j+1}-p_{j}||} = \mu_j n_T(p_j),\quad \mu_j \neq 0,\eeqq
which arises from Lagrange's multiplier theorem as within the converse implication of the proof of Proposition \ref{Prop:System}.
\epf

The next two propositions relate the positions of the hyperplanes which determine the billiard reflection rule.

\bprop\label{Prop:normalvectorsspanning}
Let $T\subset\R^n$ be a billiard table, $p=(p_1,...,p_m)$ a closed billiard trajectory on $T$, and $V\subseteq\R^n$ the affine subspace such that $T\cap V$ is the inclusion minimal affine section of $T$ that contains $p$. Further, let $V_0$ be the linear subspace parallel to $V$ satisfying $\dim V_0 = \dim V$ and let $H_1,...,H_m$ be the $T$-supporting hyperplanes through $p_1,...,p_m$ which are associated to the billiard reflection rule. Then, the convex cone spanned by the outer unit vectors
\beqq n_T(p_1),...,n_T(p_m)\eeqq
is $V_0$.
\eprop

\bpf
Since $p$ is a closed billiard trajectory on $T$ and $T\cap V$ is an affine section of $T$ that contains $p$, the vectors 
\beq p_2-p_1,...,p_m-p_{m-1},p_1-p_m\label{eq:conespanning1}\eeq
are all in $V_0$. Since $T\cap V$ is the inclusion minimal affine section that contains $p$, it follows that the convex cone spanned by the vectors in \eqref{eq:conespanning1} actually is $V_0$.

Indeed, it is
\beqq (p_2-p_1)+...+(p_m-p_{m-1})+(p_1-p_m)=0.\eeqq
For
\beqq s_j:=\frac{1}{m}\quad \forall j\in\{1,...,m\}\eeqq
we then also have
\beqq s_1(p_2-p_1)+...+s_{m-1}(p_m-p_{m-1})+s_m(p_1-p_m)=0.\eeqq
Since
\beqq \sum_{j=1}^m s_j =1 \; \text{ and } \; s_j \neq 0 \; \; \forall j\in\{1,...,m\},\eeqq
it follows that $0$ lies in the relative interior of
\beqq \conv\{p_2-p_1,...,p_m-p_{m-1},p_1-p_m\}.\eeqq
But this implies that the convex cone spanned by the vectors \eqref{eq:conespanning1} is a linear subspace of $\R^n$, more precisely, the inclusion minimal linear subspace of $\R^n$ that contains the vectors \eqref{eq:conespanning1}, i.e., $V_0$.

From the proof of Proposition \ref{Prop:System} it follows that there are $n_1,...,n_m\in S^{n-1}$ (see \eqref{eq:system3} for the definition) such that
\beqq \begin{cases} p_{j+1}-p_j=\lambda_j n_j, \; \lambda_j >0, \\ n_{j+1}-n_j = -\mu_{j+1} n_T(p_{j+1}), \; \mu_{j+1}>0, \end{cases}\eeqq
holds for all $j\in\{1,...,m\}$. From
\beqq p_{j+1}-p_j=\lambda_j n_j,\; \lambda_j>0,\eeqq
for all $j\in\{1,...,m\}$ it follows that the convex cone spanned by $n_1,...,n_m$ is $V_0$. By the above argumentation, the convex cone spanned by the vectors
\beqq n_2-n_1,...,n_m-n_{m-1},n_1-n_m\eeqq
also is a linear subspace of $\R^n$ and, therefore, equal to $V_0$. Then,
\beqq n_{j+1}-n_j=-\mu_{j+1}n_T(p_{j+1}),\;\mu_{j+1}>0,\eeqq
for all $j\in\{1,...,m\}$ implies that the convex cone spanned by the normal vectors
\beqq n_T(p_1),...,n_T(p_m)\eeqq
is $V_0$.
\epf

\bprop\label{Prop:noteinproof}
Let $T\subset\R^n$ be a billiard table, $p=(p_1,...,p_m)$ a closed billiard trajectory on $T$, and $V\subseteq\R^n$ the affine subspace such that $T\cap V$ is the inclusion minimal affine section of $T$ that contains $p$. Let $H_1,...,H_m$ be the $T$-supporting hyperplanes through $p_1,...,p_m$ which are related to the billiard reflection rule and let $H_1^+,...,H_m^+$ be the half-spaces of $\R^n$ that contain $T$ and which are bounded by $H_1,...,H_m$. Further, let $W$ be the orthogonal complement to $V$ in $\R^n$. Then, we can write
\beqq H_j=\left(H_j\cap V\right)\oplus W\;\text{ and }\; H_j^+=\left(H_j^+\cap V\right)\oplus W\eeqq
for all $j\in\{1,...,m\}$ and have that
\beqq \bigcap_{j=1}^m \left(H_j^+\cap V\right) \text{ is bounded in }V,\quad \bigcap_{j=1}^m H_j^+ \text{ is nearly bounded in }\R^n. \eeqq
\eprop

\bpf
By Proposition \ref{Prop:normalvectorsspanning}, the convex cone spanned by the outer unit vectors normal to $H_1,...,H_m$ is the linear subspace $V_0$ that underlies the affine subspace $V$. This implies, on the one hand, that we can write
\beq H_j=\left(H_j\cap V\right)\oplus W\;\text{ and }\; H_j^+=\left(H_j^+\cap V\right)\oplus W\label{eq:noteinproof0}\eeq
for all $j\in\{1,...,m\}$ and, on the other hand, that
\beq \bigcap_{j=1}^m (H_j^+\cap V)\label{eq:noteinproof1}\eeq
is bounded in $V$. The latter fact implies that there are parallel hyperplanes $H$ and $H+d$, $d\in V_0$, in $V$ such that \eqref{eq:noteinproof1} lies in between. With \eqref{eq:noteinproof0} this implies that
\beqq \bigcap_{j=1}^m H_j^+=\bigcap_{j=1}^m \left(\left(H_j^+\cap V\right)\oplus W\right) = \left(\bigcap_{j=1}^m \left(H_j^+\cap V\right)\right)\oplus W\eeqq
lies between the parallel hyperplanes
\beqq H\oplus W \; \text{ and } \; (H+d)\oplus W \eeqq
in $\R^n$ and therefore is nearly bounded in $\R^n$.
\epf

Considering Lemma \ref{Lem:Bezdek}(ii), we note that with Proposition \ref{Prop:noteinproof} we have subsequently provided a proof of Lemma \ref{Lem:Bezdek}(iii).

The following proposition is a preparation for the proof of Theorem \ref{Thm:RegularityResult1}.

\bprop\label{Prop:GenBezdek}
Let $T\subset\R^n$ be a billiard table, $p=(p_1,...,p_m)$ a closed billiard trajectory on $T$, and $V\subseteq\R^n$ the affine subspace such that $T\cap V$ is the inclusion minimal affine section of $T$ that contains $p$. Then, there is a selection
\beqq \{i_1,...,i_{\dim V+1}\}\subseteq\{1,...,m\}\eeqq
such that
\beqq \left\{p_{i_1},...,p_{i_{\dim V+1}}\right\}\in F(T).\eeqq
\eprop

\bpf
If $\dim V=n$, then the claim follows immediately by Lemma \ref{Lem:Bezdek}(i)$\&$(iii). This is also the case when $m=\dim V+1$.

\begin{figure}[h!]
\centering
\def\svgwidth{250pt}
\begingroup%
  \makeatletter%
  \providecommand\color[2][]{%
    \errmessage{(Inkscape) Color is used for the text in Inkscape, but the package 'color.sty' is not loaded}%
    \renewcommand\color[2][]{}%
  }%
  \providecommand\transparent[1]{%
    \errmessage{(Inkscape) Transparency is used (non-zero) for the text in Inkscape, but the package 'transparent.sty' is not loaded}%
    \renewcommand\transparent[1]{}%
  }%
  \providecommand\rotatebox[2]{#2}%
  \newcommand*\fsize{\dimexpr\f@size pt\relax}%
  \newcommand*\lineheight[1]{\fontsize{\fsize}{#1\fsize}\selectfont}%
  \ifx\svgwidth\undefined%
    \setlength{\unitlength}{340.67058538bp}%
    \ifx\svgscale\undefined%
      \relax%
    \else%
      \setlength{\unitlength}{\unitlength * \real{\svgscale}}%
    \fi%
  \else%
    \setlength{\unitlength}{\svgwidth}%
  \fi%
  \global\let\svgwidth\undefined%
  \global\let\svgscale\undefined%
  \makeatother%
  \begin{picture}(1,0.78120034)%
    \lineheight{1}%
    \setlength\tabcolsep{0pt}%
    \put(0,0){\includegraphics[width=\unitlength,page=1]{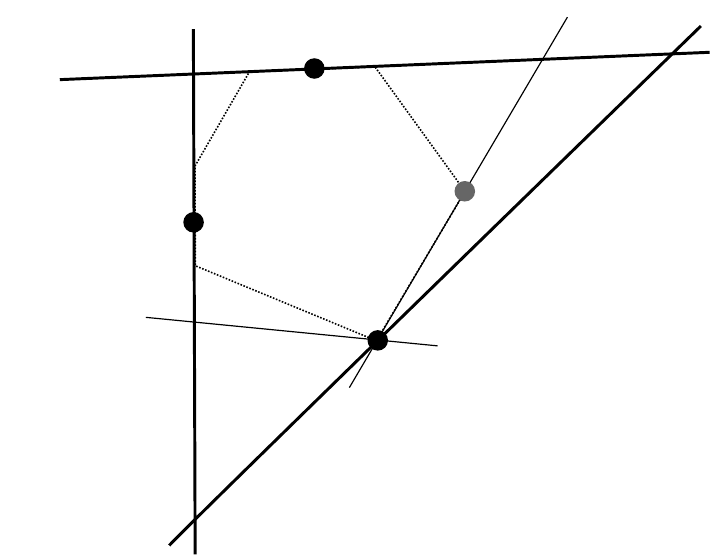}}%
    \put(0.39129674,0.71462416){\color[rgb]{0,0,0}\makebox(0,0)[lt]{\lineheight{1.25}\smash{\begin{tabular}[t]{l}$p_1=p_{i_1}$\end{tabular}}}}%
    \put(0.13122558,0.46245419){\color[rgb]{0,0,0}\makebox(0,0)[lt]{\lineheight{1.25}\smash{\begin{tabular}[t]{l}$p_2=p_{i_2}$\end{tabular}}}}%
    \put(0.52946425,0.2588915){\color[rgb]{0,0,0}\makebox(0,0)[lt]{\lineheight{1.25}\smash{\begin{tabular}[t]{l}$p_3=p_{i_3}$\end{tabular}}}}%
    \put(0.59217038,0.50619429){\color[rgb]{0,0,0}\makebox(0,0)[lt]{\lineheight{1.25}\smash{\begin{tabular}[t]{l}$p_4$\end{tabular}}}}%
    \put(-0.0012183,0.62381034){\color[rgb]{0,0,0}\makebox(0,0)[lt]{\lineheight{1.25}\smash{\begin{tabular}[t]{l}$H_1\cap V=\widetilde{H}_1$\end{tabular}}}}%
    \put(0.13022475,0.75530098){\color[rgb]{0,0,0}\makebox(0,0)[lt]{\lineheight{1.25}\smash{\begin{tabular}[t]{l}$H_2\cap V=\widetilde{H}_2$\end{tabular}}}}%
    \put(0.13441876,0.29843616){\color[rgb]{0,0,0}\makebox(0,0)[lt]{\lineheight{1.25}\smash{\begin{tabular}[t]{l}$H_3\cap V$\end{tabular}}}}%
    \put(0.75488156,0.77004931){\color[rgb]{0,0,0}\makebox(0,0)[lt]{\lineheight{1.25}\smash{\begin{tabular}[t]{l}$H_4\cap V$\end{tabular}}}}%
    \put(0.38718079,0.11937277){\color[rgb]{0,0,0}\makebox(0,0)[lt]{\lineheight{1.25}\smash{\begin{tabular}[t]{l}$\widetilde{H}_{3}$\end{tabular}}}}%
    \put(0.3680655,0.543614){\color[rgb]{0,0,0}\makebox(0,0)[lt]{\lineheight{1.25}\smash{\begin{tabular}[t]{l}$T\cap V$\end{tabular}}}}%
  \end{picture}%
\endgroup%
\caption[The choice of hyperplanes $\widetilde{H}_j$.]{Illustration of the selection of $\{p_{i_1},p_{i_2},p_{i_3}\}$ out of $\{p_1,p_2,p_3,p_4\}$ and the choice of $\bigcap_{j=1}^m\left(H_j^+\cap V\right)$-supporting hyperplanes $\widetilde{H}_1,\widetilde{H}_2,\widetilde{H}_3$ in $V$ such that $\bigcap_{j=1}^3\widetilde{H}^+_j$ is nearly bounded in $V$.}
\label{img:hyperplanechoiceeuclid}
\end{figure}

Let
\beqq \dim V\leq \min\{n-1,m-2\}.\eeqq
Since $p$ is a closed billiard trajectory on $T$, there are $T$-supporting hyperplanes $H_1,...,H_m$ through $p_1,...,p_m$ for which the billiard reflection rule is satisfied. Proposition \ref{Prop:noteinproof} implies, on the one hand, that we can write
\beqq H_j=\left(H_j\cap V\right)\oplus W\;\text{ and }\; H_j^+=\left(H_j^+\cap V\right)\oplus W\eeqq
for all $j\in\{1,...,m\}$, where $W$ is the orthogonal complement to $V$ and
\beqq H_1^+,...,H_m^+\eeqq
are the closed half-spaces of $\R^n$ defined by
\beqq \partial H_j^+=H_j \; \text{ and } \; T\subset H_j^+\eeqq
for all $j\in\{1,...,m\}$, and, on the other hand, that
\beqq \bigcap_{j=1}^m \left(H_j^+\cap V\right)\eeqq
is bounded in $V$. By Lemma \ref{Lem:Bezdek}(ii), this implies that
\beqq p\in F\left(\bigcap_{j=1}^m \left(H_j^+\cap V\right)\right).\eeqq
Then, by Lemma \ref{Lem:Bezdek}(i), there is a selection
\beqq  \{i_1,...,i_{\dim V+1}\}\subset \{1,...,m\}\eeqq
such that
\beqq \{p_{i_1},...,p_{i_{\dim V +1}}\}\in F\left(\bigcap_{j=1}^m \left(H_j^+\cap V\right)\right).\eeqq
Referring again to Lemma \ref{Lem:Bezdek}(ii), there are $\bigcap_{j=1}^m \left(H_j^+\cap V\right)$-supporting hyperplanes\footnote{It does not necessarily have to be $\widetilde{H}_j=H_{i_j}\cap V$. Further, note that due to this special situation it is not necessary to proceed to an even finer selection of $\{i_1,...,i_{\dim V +1}\}$ (and even if that were necessary, it would not impede the following argument).}
\beqq \widetilde{H}_1,...,\widetilde{H}_{\dim V +1}\eeqq
in $V$ through
\beqq p_{i_1},...,p_{i_{\dim V +1}}\eeqq
such that
\beqq \bigcap_{j=1}^{\dim V +1} \widetilde{H}_j^+\eeqq
is nearly bounded in $V$ with
\beqq \bigcap_{j=1}^m \left(H_j^+\cap V\right) \subseteq \bigcap_{j=1}^{\dim V +1} \widetilde{H}_j^+,\eeqq
where $\widetilde{H}_j^+$ is the half-space within $V$ which is bounded by $\widetilde{H}_j$ and contains $T\cap V$ for all $j\in\{1,...,\dim V +1\}$ (see Figure \ref{img:hyperplanechoiceeuclid}).
Then, this implies that
\beqq \bigcap_{j=1}^{\dim V +1}\left(\widetilde{H}_j^+ \oplus W\right)\eeqq
is nearly bounded in $\R^n$ with
\beqq T \subseteq \bigcap_{j=1}^m H_j^+ \subseteq \bigcap_{j=1}^{\dim V +1}\left(\widetilde{H}_j^+ \oplus W\right)\eeqq
and
\beqq p_{i_j}\in \widetilde{H}_j \oplus W\quad \forall j\in\{1,...,\dim V +1\}.\eeqq
By using Lemma \ref{Lem:Bezdek}(ii), this yields
\beqq \{p_{i_1},...,p_{i_{\dim V +1}}\} \in F(T).\eeqq
\epf

\section{Proof of Theorem \ref{Thm:RegularityResult1}}\label{Sec:Proof1}

In order to prove Theorem \ref{Thm:RegularityResult1}, it will be useful to formulate the following lemma:

\blem\label{Lem:boundednessopencondition}
Let $H_1^+,...,H_k^+$ be half-spaces of $\R^{d\geq 2}$ such that
\beqq H_1^+ \cap ... \cap H_k^+\eeqq
is bounded in $\R^d$ (therefore, $k\geq d+1$). Let $n_1,...,n_k$ be the outer (with respect to $H_1^+,...,H_k^+$) unit vectors normal to
\beqq H_1=\partial H_1^+,...,H_k=\partial H_k^+.\eeqq
The following holds for every $j\in\{1,...,k\}$: There is an $\eps_j>0$ such that
\beqq H_j^{pert,+}\cap \left(\bigcap_{i=1, i\neq j}^k H_i^+\right)\eeqq
is bounded in $\R^d$ for all
\beqq H_j^{pert}=\partial H_j^{pert,+}\eeqq
whose outer unit normal vector is an element of
\beqq S^{d-1} \cap B_{\eps_j}^d(n_j),\eeqq
where by $B_{\eps_j}^d(n_j)$ we denote the $d$-dimensional ball of radius $\eps_j$ centred at $n_j$.
\elem

\bpf
The statement is equivalent to the following one: Let
\beqq n_1,...,n_k\in S^{d-1}\eeqq
be unit vectors with $0$ in the interior of the convex hull
\beq \conv\{n_1,...,n_k\}.\label{eq:Lemmaconvexhull}\eeq
Then, for every $j\in\{1,...,k\}$ there is an $\eps_j>0$ such that $0$ is in the interior of
\beqq \conv\left\{n_1,...,n_{j-1},n_j^{pert},n_{j+1},...,n_k\right\}\eeqq
for every
\beqq n_j^{pert}\in S^{d-1}\cap B_{\eps_j}(n_j).\eeqq
But this is clear since for every $j\in\{1,...,k\}$ the fact
\beqq \text{\glqq $0$ is in the interior of $\conv\{n_1,...,n_k\}$\grqq}\eeqq
is invariant under small perturbations of $n_j$.

Indeed, without loss of generality we can assume $j=1$. If $0$ is in the interior of \eqref{eq:Lemmaconvexhull}, then there are $s_1,...,s_k > 0$ such that
\beqq \sum_{j=1}^k s_j = 1 \; \text{ and } \; \sum_{j=1}^k s_j n_j.\eeqq
Let
\beqq n_1^{pert}=n_1 + m \in S^{d-1}\cap B_{\eps_1}(n_1)\eeqq
be a small perturbation of $n_1$. Then, since $0$ is in the interior of \eqref{eq:Lemmaconvexhull}, we can find $\alpha_2,...,\alpha_k \in\R$ such that
\beqq -s_1 m = \sum_{j=2}^k \alpha_j n_j.\eeqq
We define
\beqq \widetilde{s}_1:=\frac{s_1}{1+\sum_{j=2}^k \alpha_j}, \quad \widetilde{s}_j:= \frac{s_j+\alpha_j}{1+\sum_{j=2}^k \alpha_j}\; \forall j\in\{2,...,k\}.\eeqq 
Then,
\beqq \sum_{j=1}^k \widetilde{s}_j = \frac{1}{1+\sum_{j=2}^k \alpha_j}\left( s_1 + \sum_{j=2}^k s_j +\sum_{j=2}^k \alpha_j\right) = 1\eeqq
and
\begin{align*}
\widetilde{s}_j n_1^{pert} + \sum_{j=2}^k \widetilde{s}_j n_j &= \frac{1}{1+\sum_{j=2}^k \alpha_j}\left(s_1 (n_1+m)+\sum_{j=2}^m (s_j + \alpha_j) n_j\right)\\
&=\frac{1}{1+\sum_{j=2}^k \alpha_j}\left(\sum_{j=1}^k s_j n_j + s_1 m + \sum_{j=2}^k \alpha_j n_j\right)\\
&=0,
\end{align*}
where $\widetilde{s}_j >0$ for all $j\in\{1,...,k\}$ since $\eps_1$ can be chosen so small that $|\alpha_1|,...,|\alpha_k|$ are sufficiently small. This implies that $0$ is in the interior of the convex hull
\beqq \conv\left\{n_1^{pert},n_2,...,n_k\right\}.\eeqq
\epf

For the proof of Theorem \ref{Thm:RegularityResult1} we finally need the following lemma.

For that, for a point $q\in \partial (T\cap V)$, where $T$ is a convex body in $\R^n$, $V$ a section of $T$, and $V_0$ the linear subspace parallel to $V$ satisfying $\dim V_0 = \dim V$, we define
\beqq N_{T\cap V}^{V_0}(q)=\left\{n\in V_0 : \langle n,y-q \rangle \leq 0 \text{ for all }y\in T\cap V\right\}.\eeqq

\blem\label{Lem:normalconerelation}
Let $T\subset\R^n$ be a convex body and $\{p_1,...,p_m\}$ a set of boundary points of $T$ such that the outer normal unit vectors
\beqq n_T(p_1),...,n_T(p_m) \; \text{ in } \; N_T(p_1),...,N_T(p_m)\eeqq
lie in $V_0$ by which we denote the linear subspace parallel to the inclusion minimal affine section $V$ of $T$ that contains $\{p_1,...,p_m\}$ satisfying $\dim V_0 = \dim V$. Then, we have
\beqq N_T(p_j)\cap V_0 = N_T(p_j) \cap N_{T\cap V}^{V_0}(p_j)\eeqq
for all $j\in\{1,...,m\}$.
\elem

\bpf
From
\beqq N_{T\cap V}^{V_0}(p_j) \subseteq V_0\eeqq
for all $j\in\{1,...,m\}$ it follows
\beqq N_T(p_j)\cap V_0 \supseteq N_T(p_j) \cap N_{T\cap V}^{V_0}(p_j)\eeqq
for all $j\in\{1,...,m\}$.

Let $j\in\{1,...,m\}$ be arbitrarily chosen. Let $n$ be a nonzero vector in
\beqq N_T(p_j)\cap V_0.\eeqq
Then,
\beqq n\in N_T(p_j),\text{ i.e., }\; \langle n,x-p_j\rangle \leq 0 \;\; \forall x\in T,\eeqq
and $n\in V_0$. Because of
\beqq T\cap V\subseteq T,\eeqq
this implies
\beqq \langle n,x-p_j \rangle \leq 0 \;\; \forall x\in T\cap V,\; n\in V_0,\; n\in N_T(p_j).\eeqq
From that, we conclude
\beqq n\in N_{T\cap V}^{V_0}(p_j) \; \text{ and } \; n\in N_T(p_j)\eeqq
and, therefore,
\beqq n\in N_{T\cap V}^{V_0}(p_j)\cap N_T(p_j).\eeqq
Consequently,
\beqq N_T(p_j)\cap V_0 \subseteq N_T(p_j) \cap N_{T\cap V}^{V^0}(p_j)\eeqq
for all $j\in\{1,...,m\}$.
\epf

We come to the proof of Theorem \ref{Thm:RegularityResult1}:

\bpf[of Theorem \ref{Thm:RegularityResult1}]
Let $p=(p_1,...,p_m)$ be a length-minimizing closed billiard trajectory on $T$ and $V\subseteq\R^n$ the affine subspace such that $T\cap V$ is the inclusion minimal affine section that contains $p$, i.e.,
\beqq \dim V \leq \min\{n,m-1\}.\eeqq

Then, by Proposition \ref{Prop:GenBezdek}, there is a selection
\beqq \{i_1,...,i_{\dim V+1}\}\subseteq \{1,...,m\}\eeqq
with
\beqq \left\{p_{i_1},...,p_{i_{\dim V+1}}\right\}\in F(T).\eeqq
Without loss of generality, we can assume
\beqq i_1 <...< i_{\dim V+1}\eeqq
and define the closed polygonal curve
\beqq \widetilde{p}:=\left(p_{i_1},...,p_{i_{\dim V+1}}\right).\eeqq
For
\beqq \dim V +1 < m\eeqq
it follows
\beqq \ell(\widetilde{p})< \ell (p).\eeqq
But with Theorem \ref{Thm:Bezdek}, this is a contradiction to the minimality of $p$. Therefore, it follows
\beqq \dim V=m-1.\eeqq

Let us denote by $H_1^+,...,H_m^+$ the closed half-spaces of $\R^n$ defined by
\beqq \partial H_j^+=H_j \; \text{ and} \; T\subseteq H_j^+\eeqq
for all $j\in\{1,...,m\}$, where by $H_1,...,H_m$ we denote the $T$-supporting hyperplanes through $p_1,...,p_m$ which are related to the billiard reflection rule in these points. By Proposition \ref{Prop:noteinproof}, we conclude that we can write
\beq H_j=(H_j\cap V)\oplus W\; \text{ and }\;H_j^+=\left(H_j^+\cap V\right)\oplus W\label{eq:regularsplitting}\eeq
for all $j\in\{1,...,m\}$, where $W$ is the orthogonal complement to $V$ in $\R^n$, and that
\beqq \bigcap_{j=1}^m \left(H_j^+\cap V\right) \text{ is bounded in }V, \; \bigcap_{j=1}^m H_j^+ \text{ is nearly bounded in }\R^n.\eeqq

Let $n_1,...,n_m$ be the outer unit vectors normal to $H_1,...,H_m$. Then, it follows by \eqref{eq:regularsplitting} that
\beqq n_j\in N_T(p_j)\cap V_0\eeqq
and, therefore,
\beqq \dim (N_T(p_j)\cap V_0)\geq 1\eeqq
for all $j\in\{1,...,m\}$.

Let us assume there is an $i\in\{1,...,m\}$ such that
\beqq \dim (N_T(p_i)\cap V_0)>1.\eeqq
Then, using Lemma \ref{Lem:normalconerelation}, i.e.,
\beqq N_T(p_i)\cap V_0 = N_{T}(p_i)\cap N_{T\cap V}^{V_0}(p_i),\eeqq
it follows
\beqq \dim \left(N_T(p_i)\cap N_{T\cap V}^{V_0}(p_i)\right)>1,\eeqq
and, because of Lemma \ref{Lem:boundednessopencondition} (for $d=m-1$ and $k=m$), we can find a unit vector
\beq n_i^{pert}\in N_{T}(p_i)\cap N_{T\cap V}^{V_0}(p_i)\; \text{ with } \; n_i^{pert}\neq n_i \label{eq:choicepiperturbed}\eeq
such that 
\beq H_{i,V}^{pert,+}\cap \left(\bigcap_{j=1,j\neq i}^m\left(H_{j}^+\cap V\right)\right)\label{eq:regularityProp2}\eeq
remains bounded in $V$, where we denote by $H_{i,V}^{pert,+}$ the closed half-space of $V$ that contains $T\cap V$ and which is bounded by $H_{i,V}^{pert}$ which is the hyperplane in $V$ through $p_i$ that is normal to $n_i^{pert}$. Since by Proposition \ref{Prop:onlyonehyperplane} the billiard reflection rule in $p_i$ (as bouncing point of the closed billiard trajectory $p$ in $T\cap V$; see Proposition \ref{Prop:SectionInvariance}) is no longer satisfied with respect to the perturbed hyperplane $H_{i,V}^{pert}$, the bouncing point $p_i$ can be moved along $H_{i,V}^{pert}$, say to $p_i^*$, in order to reduce the length of the polygonal curve segment
\beqq (p_{i-1},p_i,p_{i+1}).\eeqq
We define the closed polygonal curve
\beqq \widetilde{p}:=(p_1,...,p_{i_1},p_i^*,p_{i+1},...,p_m),\eeqq
for which with the boundedness of \eqref{eq:regularityProp2} in $V$ we conclude
\beqq \widetilde{p}\in F(T\cap V)\eeqq
by Lemma \ref{Lem:Bezdek}(ii). Now, we argue that $\widetilde{p}\in F(T)$: With the boundedness of \eqref{eq:regularityProp2} in $V$ it follows with
\beq H_i^{pert}:=H_{i,V}^{pert}\oplus W,\quad H_i^{pert,+}:=H_{i,V}^{pert,+}\oplus W\label{eq:regularityProp21}\eeq
and \eqref{eq:regularsplitting} the nearly boundedness of
\beq H_i^{pert,+}\cap \left(\bigcap_{j=1,j\neq i}^m H_j^+\right)\label{eq:regularityProp3}\eeq
in $\R^n$.

Indeed, when the intersection in \eqref{eq:regularityProp2} is bounded in $V$, then there is a hyperplane $H$ in $V$ such that the intersection lies between $H$ and $H+d$ for an appropriate $d\in V_0$. Then, it follows with \eqref{eq:regularsplitting} and \eqref{eq:regularityProp21} that
\begin{align*}
&H_i^{pert,+}\cap \left(\bigcap_{j=1,j\neq i}^m H_j^+\right)\\
= & \left(H_{i,V}^{pert,+}\oplus W\right)\cap \left(\bigcap_{j=1,j\neq i}^m \left(\left(H_j^+\cap V\right)\oplus W\right)\right)\\
=&\left(H_{i,V}^{pert,+}\cap \left(\bigcap_{j=1,j\neq i}^m \left(H_j^+\cap V\right)\right)\right)\oplus W
\end{align*}
lies between the hyperplanes
\beqq H\oplus W \; \text{ and } \; (H+d)\oplus W.\eeqq

Since $H_i^{pert}$ is a $T$-supporting hyperplane through $p_i$ (what follows from the fact that by \eqref{eq:choicepiperturbed} its outer unit normal vector $n_i^{pert}$ is an element of $N_T(p_i)$), we conclude that $T$ is a subset of the intersection in \eqref{eq:regularityProp3}. Then, it follows from the nearly boundedness (in $\R^n$) of the intersection in \eqref{eq:regularityProp3} together with Lemma \ref{Lem:Bezdek}(ii) that
\beqq \widetilde{p}\in F(T).\eeqq
By referring to Theorem \ref{Thm:Bezdek}, from
\beqq \ell(\widetilde{p})<\ell(p)\eeqq
we derive a contradiction to the minimality of $p$.

Therefore, it follows
\beqq \dim (N_T(p_i)\cap V_0)=1.\eeqq
\epf

\section{Examples}\label{Sec:Examples}

In what follows, we discuss three relevant examples showing:
\begin{itemize}
\item[(A)] The statement of Proposition \ref{Prop:GenBezdek} is not true when requiring $p$ just to be a closed polygonal curve in $F(T)$ (and not a closed billiard trajectory on $T$).
\item[(B)] Let $T\subset\R^n$ be a billiard table and $p=(p_1,...,p_m)$ a length-minimizing closed billiard trajectory on $T$. Let $T\cap V$ be the inclusion minimal affine section of $T$ that contains $p$. Then, $p$ may not even locally minimize the length of closed polygonal curves in $F(T\cap V)$.
\item[(C)] A length-minimizing closed billiard trajectory in $T\subset\R^n$ may not be regular within the inclusion minimal affine section of $T$ that contains the billiard trajectory. This can even appear for the unique length-minimizing closed billiard trajectory.
\end{itemize}

\noindent\underline{Ad (A)}: We consider the following example\footnote{This example arose from the hints of A. Abbondandolo.} in $\R^5$: We start from four convex bodies $K_1,K_2,K_3,K_4$ in $\R^3$ with the following two properties:
\begin{itemize}
\item[(a)] the intersection of all of them is empty;
\item[(b)] the intersection of any three of them has non-empty interior.
\end{itemize}
One has these examples because Helly's theorem is sharp. Then, we consider the four vertices of a square in $\R^2$:
\beqq  v_1=(0,0),\; v_2=(1,0),\; v_3=(1,1),\; v_4=(0,1).\eeqq
Now, let $T$ be the convex hull in
\beqq \R^5=\R^2\times \R^3\eeqq
of the union of the following four sets:
\beqq \{v_1\}\times K_1,\; \{v_2\}\times K_2,\; \{v_3\}\times K_3,\; \{v_4\}\times K_4.\eeqq
$T$ projects onto the square, but (a) implies that each section of $T$ that is parallel to $\R^2\times \{0\}$ has area smaller than $1$.

Indeed, we take any $w\in\R^3$ and look at the section
\beq T \cap (\R^2 \times {w}).\label{eq:ExAlberto0}\eeq
By the definition of $T$, the points in this section are of the form $(v,w)$ with
\beq v = \sum_j \lambda_j v_j\; \text{ and }\; w=  \sum_j \lambda_j w_j, \label{eq:ExAlberto}\eeq
where $w_j$ is in $K_j$ and the $\lambda_j$s are positive and add up to one. In particular, this section is contained in $Q \times {w}$, where $Q$ denotes the square with vertices $v_j$. This section cannot contain all the four points $(v_j,w)$. In fact, assume that it contains the point $(v_1,w)$. Then, \eqref{eq:ExAlberto} and the fact that $v_1$ is an extremal point of $Q$ imply that $w$ belongs to $K_1$, as all the $\lambda_j$s with $j>1$ must vanish in \eqref{eq:ExAlberto}. Since any given $w$ in $\R^3$ belongs to at most three of the $K_j$s, the claim is proven. Being a closed set that is contained in $Q \times {w}$ and does not contain $(v_j,w)$ for at least one $j$, the section \eqref{eq:ExAlberto0} has area strictly smaller than $1$. Since the area of the intersection of a convex body with $\R^2 \times {w}$ is an upper semi-continuous function of $w$, all the sections of $T$ by planes parallel to $\R^2 \times \{0\}$ have area less than $A$ for some $A<1$.

We choose as $p_1,p_2,p_3,p_4$ the points
\beqq p_1=(1-t)v_1,\;p_2=(1-t)v_2,\; p_3=(1-t)v_3,\; p_4=(1-t)v_4\eeqq
for $t>0$ so small that the area of the square with the vertices $p_1,p_2,p_3,p_4$ is larger than $A$ (while smaller than $1$). Then,
\beqq \{p_1,p_2,p_3,p_4\}\in F(T)\eeqq
because any translation of these points will enclose a square that is too big to be contained in a section of $T$ parallel to $\R^2\times \{0\}$. However, any triplet from $\{p_1,p_2,p_3,p_4\}$ is not in $F(T)$: Consider without loss of generality the triplet $\{p_1,p_2,p_3\}$. By (b), there is a point $w$ in the interior of
\beqq K_1\cap K_2 \cap K_3.\eeqq
Then, the section $\R^2\times \{w\}$ contains a translated copy of the triangle with vertices $v_1,v_2,v_3$ and $\{p_1,p_2,p_3\}$ can be translated into the interior of such triangle, and, hence, into the interior of $T$. This proves the claim.

We remark that the proof of Proposition \ref{Prop:GenBezdek} does not apply to this example due to the lack of applicability of Proposition \ref{Prop:noteinproof} which itself goes back to Proposition \ref{Prop:normalvectorsspanning} (for which being a closed billiard trajectory is essential beyond the sole being of a closed polygonal curve that cannot be translated into the interior of $T$.)

\noindent\underline{Ad (B)}: Let $T\subset\R^3$ be the convex hull of the points
\beqq (0,0,0), (1,0,0), \left(\frac{1}{2},0,\frac{\sqrt{3}}{2}\right), (0,-2,0), (1,-2,0).\eeqq
One checks that $p=(p_1,p_2,p_3)$ with
\beqq p_1=\left(\frac{1}{2},-1,0\right),\;p_2=\left(\frac{1}{4},-1,\frac{\sqrt{3}}{4}\right),\; p_3=\left(\frac{3}{4},-1,\frac{\sqrt{3}}{4}\right)\eeqq
is a length-minimizing closed billiard trajectory on $T$: It is the Fagnano triangle of the affine section $T\cap X_{1,3}$ translated by $(0,-1,0)$, where, in general, we denote by $X_{i,j}$ the $(x_i,x_j)$-plane of $\R^3_{(x_1,x_2,x_3)}$. If we define
\beqq V:=X_{1,3}+(0,-1,0),\eeqq
then $T\cap V$ is the inclusion minimal affine section of $T$ that contains $p$. One checks that
\beqq \widetilde{p}=( p_1,\widetilde{p}_2) \; \text{ with } \; \widetilde{p}_2=\left(\frac{1}{2},-1,\frac{\sqrt{3}}{4}\right)\eeqq
is a length-minimizing closed billiard trajectory in $T\cap V$ (with $\widetilde{p}\notin F(T)$) but
\beqq \ell(\widetilde{p})<\ell(p).\eeqq
Therefore, $p$ is not a length-minimizing closed billiard trajectory in $T\cap V$.

\begin{figure}[h!]
\centering
\def\svgwidth{220pt}
\begingroup%
  \makeatletter%
  \providecommand\color[2][]{%
    \errmessage{(Inkscape) Color is used for the text in Inkscape, but the package 'color.sty' is not loaded}%
    \renewcommand\color[2][]{}%
  }%
  \providecommand\transparent[1]{%
    \errmessage{(Inkscape) Transparency is used (non-zero) for the text in Inkscape, but the package 'transparent.sty' is not loaded}%
    \renewcommand\transparent[1]{}%
  }%
  \providecommand\rotatebox[2]{#2}%
  \newcommand*\fsize{\dimexpr\f@size pt\relax}%
  \newcommand*\lineheight[1]{\fontsize{\fsize}{#1\fsize}\selectfont}%
  \ifx\svgwidth\undefined%
    \setlength{\unitlength}{233.44887582bp}%
    \ifx\svgscale\undefined%
      \relax%
    \else%
      \setlength{\unitlength}{\unitlength * \real{\svgscale}}%
    \fi%
  \else%
    \setlength{\unitlength}{\svgwidth}%
  \fi%
  \global\let\svgwidth\undefined%
  \global\let\svgscale\undefined%
  \makeatother%
  \begin{picture}(1,0.62739337)%
    \lineheight{1}%
    \setlength\tabcolsep{0pt}%
    \put(0,0){\includegraphics[width=\unitlength,page=1]{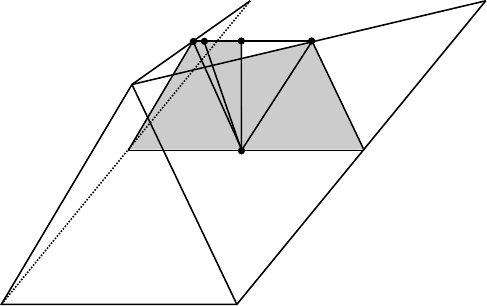}}%
    \put(0.49185662,0.28813436){\color[rgb]{0,0,0}\makebox(0,0)[lt]{\lineheight{1.25}\smash{\begin{tabular}[t]{l}$p_1$\end{tabular}}}}%
    \put(0.34933445,0.56044338){\color[rgb]{0,0,0}\makebox(0,0)[lt]{\lineheight{1.25}\smash{\begin{tabular}[t]{l}$p_2$\end{tabular}}}}%
    \put(0.63224076,0.57009773){\color[rgb]{0,0,0}\makebox(0,0)[lt]{\lineheight{1.25}\smash{\begin{tabular}[t]{l}$p_3$\end{tabular}}}}%
    \put(0.48996445,0.5642216){\color[rgb]{0,0,0}\makebox(0,0)[lt]{\lineheight{1.25}\smash{\begin{tabular}[t]{l}$\widetilde{p}_2$\end{tabular}}}}%
    \put(0.40136932,0.5836285){\color[rgb]{0,0,0}\makebox(0,0)[lt]{\lineheight{1.25}\smash{\begin{tabular}[t]{l}$\widebar{p}_2^\delta$\end{tabular}}}}%
    \put(0,0){\includegraphics[width=\unitlength,page=2]{Prop2.pdf}}%
    \put(0.56904271,0.36237441){\color[rgb]{0,0,0}\makebox(0,0)[lt]{\lineheight{1.25}\smash{\begin{tabular}[t]{l}$T\cap V$\end{tabular}}}}%
    \put(0,0){\includegraphics[width=\unitlength,page=3]{Prop2.pdf}}%
  \end{picture}%
\endgroup%
\caption{The closed billiard trajectory $p=(p_1,p_2,p_3)$ is a length-minimizing billiard trajectory in $T$, but not in $T\cap V$ (the closed billiard trajectory $(p_1,\widetilde{p}_2)$ is shorter than $p$). Furthermore, the closed polygonal line $(p_1,\widebar{p}_2^\delta,p_3)$ is shorter than $p$ for every $\delta > 0$ while it approximates $p$ for $\delta\rightarrow 0$.}
\label{img:Proposition2}
\end{figure}

The closed billiard trajectory $p$ may not even locally minimize the length of closed polygonal curves in $F(T\cap V)$. To see this, we slightly move $p_2$ clockwise along $\partial (T\cap V)$. We denote this slightly perturbed $p_2$ by
\beqq \widebar{p}^\delta_2=\left(\frac{1}{4}+\delta,-1,\frac{\sqrt{3}}{4}\right), \; \delta >0\text{ small.}\eeqq
The closed polygonal curve
\beqq \widebar{p}^\delta=\left(\widebar{p}_1,\widebar{p}^\delta_2,\widebar{p}_3\right) \; \text{ with } \; \widebar{p}_1=p_1 \; \text{ and } \;\widebar{p}_3=p_3\eeqq
satisfies
\beqq \widebar{p}^\delta\in F(T\cap V),\;\; \widebar{p}^\delta\notin F(T),\;\;\ell(\widebar{p}^\delta)<\ell(p)\; \text{(for small $\delta >0$)}\eeqq
and $\widebar{p}^\delta$ converges with respect to the Hausdorff distance to $p$ for $\delta \rightarrow 0$.

\noindent\underline{Ad (C)}: Let $T_\eps\subset\R^3$, $\eps >0$ small, be the convex polytope given by the vertices
\begin{align*}
&(0,0,0), \left(-\frac{1}{2}+\frac{\eps}{\sqrt{3}},0,\frac{\sqrt{3}}{2}-\eps\right),\left(-\frac{1}{2}+\frac{\eps}{\sqrt{3}},0,\frac{\sqrt{3}}{2}\right), \left(\frac{1}{2}-\frac{\eps}{\sqrt{3}},0,\frac{\sqrt{3}}{2}-\eps\right),\\
&(0,-2,0), \left(-\frac{1}{2}+\frac{\eps}{\sqrt{3}},-2,\frac{\sqrt{3}}{2}-\eps\right),\left(\frac{1}{2}-\frac{\eps}{\sqrt{3}},-2,\frac{\sqrt{3}}{2}-\eps\right),\left(\frac{1}{2}-\frac{\eps}{\sqrt{3}},-2,\frac{\sqrt{3}}{2}\right).
\end{align*}

\begin{figure}[h!]
\centering
\def\svgwidth{300pt}
\begingroup%
  \makeatletter%
  \providecommand\color[2][]{%
    \errmessage{(Inkscape) Color is used for the text in Inkscape, but the package 'color.sty' is not loaded}%
    \renewcommand\color[2][]{}%
  }%
  \providecommand\transparent[1]{%
    \errmessage{(Inkscape) Transparency is used (non-zero) for the text in Inkscape, but the package 'transparent.sty' is not loaded}%
    \renewcommand\transparent[1]{}%
  }%
  \providecommand\rotatebox[2]{#2}%
  \newcommand*\fsize{\dimexpr\f@size pt\relax}%
  \newcommand*\lineheight[1]{\fontsize{\fsize}{#1\fsize}\selectfont}%
  \ifx\svgwidth\undefined%
    \setlength{\unitlength}{342.18559379bp}%
    \ifx\svgscale\undefined%
      \relax%
    \else%
      \setlength{\unitlength}{\unitlength * \real{\svgscale}}%
    \fi%
  \else%
    \setlength{\unitlength}{\svgwidth}%
  \fi%
  \global\let\svgwidth\undefined%
  \global\let\svgscale\undefined%
  \makeatother%
  \begin{picture}(1,0.9140936)%
    \lineheight{1}%
    \setlength\tabcolsep{0pt}%
    \put(0,0){\includegraphics[width=\unitlength,page=1]{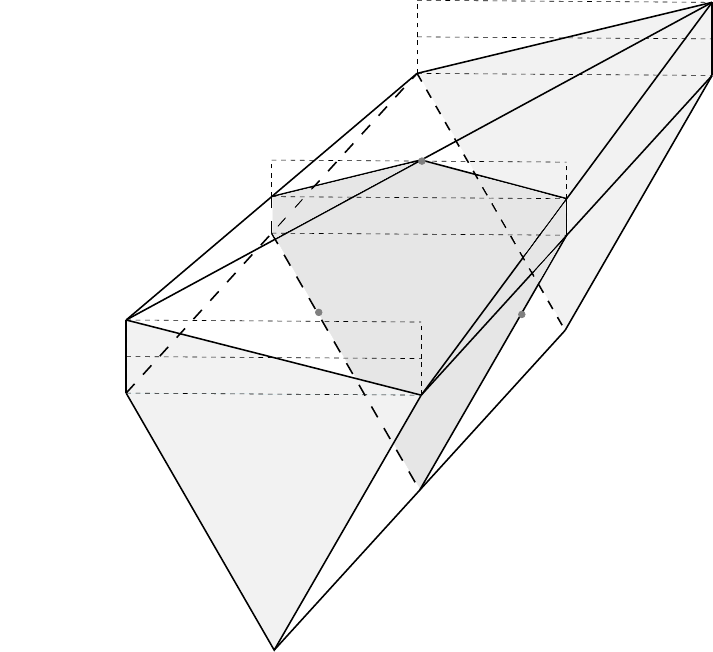}}%
    \put(0.5258885,0.52722067){\color[rgb]{0,0,0}\makebox(0,0)[lt]{\lineheight{1.25}\smash{\begin{tabular}[t]{l}$T_\eps\cap V$\end{tabular}}}}%
    \put(0.56840172,0.71333412){\color[rgb]{0,0,0}\makebox(0,0)[lt]{\lineheight{1.25}\smash{\begin{tabular}[t]{l}$p_3$\end{tabular}}}}%
    \put(0,0){\includegraphics[width=\unitlength,page=2]{Counterexamplenew_grey.pdf}}%
    \put(0.04196122,0.45999153){\color[rgb]{0,0,0}\makebox(0,0)[lt]{\lineheight{1.25}\smash{\begin{tabular}[t]{l}$\frac{\sqrt{3}}{2}$\end{tabular}}}}%
    \put(-0.00121291,0.41048298){\color[rgb]{0,0,0}\makebox(0,0)[lt]{\lineheight{1.25}\smash{\begin{tabular}[t]{l}$\frac{\sqrt{3}}{2}-\frac{\eps}{2}$\end{tabular}}}}%
    \put(0.00034539,0.36348568){\color[rgb]{0,0,0}\makebox(0,0)[lt]{\lineheight{1.25}\smash{\begin{tabular}[t]{l}$\frac{\sqrt{3}}{2}-\eps$\end{tabular}}}}%
    \put(0,0){\includegraphics[width=\unitlength,page=3]{Counterexamplenew_grey.pdf}}%
  \end{picture}%
\endgroup%
\caption{For $\eps >0$ small enough, the closed billiard trajectory $p=(p_1,p_2,p_3)$ is the unique length-minimizing billiard trajectory in $T_\eps$. It has the property that it is not regular within $T_\eps \cap V$ which is the inclusion minimal affine section of $T_\eps$ that contains $p$.}
\label{img:Counterexample}
\end{figure}

We claim that for sufficiently small $\eps >0$, the unique length-minimizing closed billiard trajectory in $T_\eps$ is given by $p=( p_1,p_2,p_3)$ with
\beqq p_1=\left(-\frac{1}{4},-1,\frac{\sqrt{3}}{4}\right),\; p_2=\left(\frac{1}{4},-1,\frac{\sqrt{3}}{4}\right),\; p_3=\left(0,-1,\frac{\sqrt{3}}{2}\right).\eeqq
Moreover, we claim that $p$ is not regular within the inclusion minimal affine section of $T_\eps$ that contains $p$.

Indeed, $p$ is contained in the affine section $T_\eps\cap V$ with
\beqq V:=X_{1,3}+(0,-1,0).\eeqq
$T_\eps\cap V$ is given by the vertices
\begin{align*}
&(0,-1,0),\left(-\frac{1}{2}+\frac{\eps}{\sqrt{3}},-1,\frac{\sqrt{3}}{2}-\eps\right), \left(-\frac{1}{2}+\frac{\eps}{\sqrt{3}},-1,\frac{\sqrt{3}}{2}-\frac{\eps}{2}\right),\\
&\left(0,-1,\frac{\sqrt{3}}{2}\right),\left(\frac{1}{2}-\frac{\eps}{\sqrt{3}},-1,\frac{\sqrt{3}}{2}-\eps\right),\left(\frac{1}{2}-\frac{\eps}{\sqrt{3}},-1,\frac{\sqrt{3}}{2}-\frac{\eps}{2}\right)
\end{align*}
and is subset of the equilateral triangle $\Delta_1$ given by the vertices
\beqq (0,-1,0),\left(-\frac{1}{2},-1,\frac{\sqrt{3}}{2}\right), \left(\frac{1}{2},-1,\frac{\sqrt{3}}{2}\right).\eeqq
The closed billiard trajectory $p$ is coinciding with the Fagnano triangle of $\Delta_1$ which is the unique length-minimizing closed billiard trajectory in $\Delta_1$ (with length $\frac{3}{2}$). By construction, $p$ is a closed billiard trajectory in $T_\eps$ as well as in $T_\eps\cap V$: The hyperplanes in $\R^3$, respectively in $V$, related to the billiard reflection rule are the one which are normal to the bisectors of the polygonal curve segments
\beqq (p_1,p_2,p_3),\;\; (p_2,p_3,p_1),\;\text{ and }\; (p_3,p_1,p_2).\eeqq

Let $F_1$ be the facet of $T_\eps$ given by the vertices
\beqq \left(-\frac{1}{2}+\frac{\eps}{\sqrt{3}},0,\frac{\sqrt{3}}{2}\right), \left(\frac{1}{2}-\frac{\eps}{\sqrt{3}},0,\frac{\sqrt{3}}{2}-\eps\right),\left(\frac{1}{2}-\frac{\eps}{\sqrt{3}},-2,\frac{\sqrt{3}}{2}\right),\eeqq
and $F_2$ the facet of $T_\eps$ given by the vertices
\beqq \left(-\frac{1}{2}+\frac{\eps}{\sqrt{3}},0,\frac{\sqrt{3}}{2}\right), \left(\frac{1}{2}-\frac{\eps}{\sqrt{3}},-2,\frac{\sqrt{3}}{2}\right), \left(-\frac{1}{2}+\frac{\eps}{\sqrt{3}},-2,\frac{\sqrt{3}}{2}-\eps\right).\eeqq

Then, $p$ uniquely minimizes the length over all closed billiard trajectories in $T_\eps$ which are contained in the affine sections
\beqq T_\eps \cap V^a,\quad V^a:=X_{1,2}+(0,-a,0),\quad a\in[0,2],\eeqq
(we note: $V^1=V$) as consequence of the tilt of $F_1$ and $F_2$. Because of this, all bouncing points of closed billiard trajectories in $T_\eps\cap V^a$ have to be on the boundary of the equilateral triangle $\Delta_1^a$ given by the vertices
\beqq (0,-a,0),\left(-\frac{1}{2},-a,\frac{\sqrt{3}}{2}\right),\left(\frac{1}{2},-a,\frac{\sqrt{3}}{2}\right)\eeqq
(we note: $\Delta_1^{1}=\Delta_1$) and the length-minimizing closed billiard trajectories in $\Delta_1^a$ are its Fagnano triangles. However, we have to be careful since there are closed billiard trajectories on $T_\eps$ which lie within the affine sections $T_\eps\cap V^a$ with bouncing points on $\Delta_1^a$, but which have the property that they are in
\beqq F\left(T_\eps\cap V^a\right)\setminus F\left(\Delta_1^a\right).\eeqq
This is for instance the case for the closed billiard trajectories given by the bouncing points
\beqq \left(-\frac{1}{2}+\frac{\eps}{\sqrt{3}},-a,\frac{\sqrt{3}}{2}-\eps\right)\; \text{ and }\; \left(\frac{1}{2}-\frac{\eps}{\sqrt{3}},-a,\frac{\sqrt{3}}{2}-\eps\right).\eeqq
However, the length of the billiard trajectories of this kind converges to
\beqq 2>\frac{3}{2}=\ell(p)\eeqq
for $\eps\rightarrow 0$.

We show that all other closed billiard trajectories in $T_\eps$ have a length greater than $\frac{3}{2}$: For that, we first notice that whenever we consider a closed billiard trajectory in $T_\eps$ with one bouncing point on the front facet
\beqq T_\eps\cap X_{1,3}\eeqq
and another one on the back facet
\beqq T_\eps\cap (X_{1,3}+(0,-2,0))\eeqq
of $T_\eps$, then it has a length greater or equal than $4$. Every other closed billiard trajectory in $T_\eps$ either has two bouncing points, one on $\mathring{F}_1$ or $\mathring{F}_2$ and the other on the back- or front facet, respectively, or it has one bouncing point on $\mathring{F}_1$ and another on $\mathring{F}_2$ while it has non-empty intersection either with the front- or with the back facet. This follows from the fact that whenever we have a closed billiard trajectory with bouncing point on $\mathring{F}_1$ ($\mathring{F}_2$) and no other on $\mathring{F}_2$ ($\mathring{F}_1$), then at least one of the previous and following bouncing points is on the front- and back facet (or the other way round), respectively. In the case of a closed billiard trajectory $p'$ with two bouncing points, one on $\mathring{F}_1$ or on $\mathring{F}_2$ and the other on the back- or the front facet, respectively, we have
\beqq \ell(p')=\min_{x\in G_1, y\in \mathring{F}_1}2|x-y|>\frac{3}{2}=\ell(p)\;\quad (\eps >0 \text{ small})\eeqq
or
\beqq \ell(p')=\min_{x\in G_2, y\in\mathring{F}_2}2|x-y|>\frac{3}{2}=\ell(p)\;\quad (\eps>0 \text{ small})\eeqq
respectively, where we denote by $G_1$ the line segment from
\beqq (0,-2,0)\;\text{ to }\; \left(-\frac{1}{2}+\frac{\eps}{\sqrt{3}},-2,\frac{\sqrt{3}}{2}-\eps\right)\eeqq
and by $G_2$ the line segment from
\beqq (0,0,0)\; \text{ to } \; \left(\frac{1}{2}-\frac{\eps}{\sqrt{3}},0,\frac{\sqrt{3}}{2}-\eps\right),\eeqq
since
\beqq \min_{x\in G_1,y\in\mathring{F}_1}|x-y| = \min_{x\in G_2,y\in\mathring{F}_2}|x-y|\rightarrow \frac{\sqrt{3}}{2}\quad (\eps\rightarrow 0).\eeqq
In the case of a closed billiard trajectory $p'$ with one bouncing point on $\mathring{F}_1$ and another on $\mathring{F}_2$ while it has non-empty intersection either with the front- or with the back facet, we have
\beqq \ell(p')\geq \min_{x\in \mathring{F}_1,y\in\mathring{F}_2,z\in G}\left(|x-z|+|z-y|\right)>\frac{3}{2}=\ell(p),\eeqq
where we denote by $G$ the line segment from
\beqq (0,0,0) \; \text{ to } \; (0,-2,0).\eeqq

The closed billiard trajectory $p$ is not regular in $T_\eps\cap V$ since the normal cone
\beqq N_{T_\eps\cap V}(p_3)\eeqq
is two-dimensional, i.e., $p_3$ is a non-smooth boundary point of $T_\eps \cap V$.

We clearly see why the argument used in the proof of Theorem \ref{Thm:RegularityResult1} does not apply to this example. The $T_\eps$-supporting hyperplane through $p_3$ for which the billiard reflection rule is satisfied is
\beqq H_3:=X_{1,2}+\left(0,0,\frac{\sqrt{3}}{2}\right).\eeqq
The only way of perturbing $H_3$ to $H_3^{pert}$ as required within the proof of Theorem \ref{Thm:RegularityResult1} is by tilting it around the axis through the points
\beqq \left(-\frac{1}{2}+\frac{\eps}{\sqrt{3}},0,\frac{\sqrt{3}}{2}\right)\; \text{ and } \;  \left( \frac{1}{2}-\frac{\eps}{\sqrt{3}},-2,\frac{\sqrt{3}}{2}\right).\eeqq
But in any case,
\beqq H_1^{+}\cap H_2^{+}\cap H_3^{pert,+}\eeqq
cannot be nearly bounded in $\R^3$ (when $H_1$ and $H_2$ denote the uniquely determined $T_\eps$-supporting hyperplanes through $p_1$ and $p_2$). Therefore, it is not possible to construct a closed polygonal curve
\beqq \widetilde{p}\in F(T_\eps\cap V)\; \text{ with } \; \ell(\widetilde{p})\leq \ell(p) \; \text{ and } \; \widetilde{p}\neq p\eeqq
while guaranteeing
\beqq \widetilde{p}\in F(T_\eps).\eeqq

We remark that $T_\eps$ is a convex polytope. Nevertheless, it can be made strictly convex without losing the relevant properties.

\section{Constructing closed billiard trajectories}\label{Sec:Construction}

\subsection{Constructing closed regular billiard trajectories}\label{Subsec:constructingregular}

Throughout this section, we let $T \subset \mathbb{R}^n$ be a convex polytope that plays the role of the billiard table. In this section, we will concentrate on closed regular billiard trajectories on $P$. Regular billiard trajectories are especially interesting given the fact that they are the classical ones, i.e., billiard trajectories with bouncing points on the billiard table cushions while not running into the billiard table pockets (of course, here we mean its higher dimensional generalizations, where we identify every $j$-face of $P$ for $j<n+1$ as a pocket and every facet of $P$ as a cushion of the billiard table $P$).

By doing this, we extend the approach of finding closed regular billiard trajectories on two dimensional convex polytopes described in \cite{AlkoumiSchlenk2014} to higher dimensions. Therein, the main tool was the property of the so called Fagnano triangles representing the uniquely determined closed regular billiard trajectories in acute triangles. Since there are no higher-dimensional analogues of the Fagnano triangles for the extension, we though have to generalize the approach: it will turn out that we can make use of the reformulation of the billiard reflection rule in Proposition \ref{Prop:System}.

We note that in higher dimensions, it is much more complicated to construct closed non-regular billiard trajectories. However, the special case of closed non-regular billiard trajectories with two bouncing points is not a challenge. For that, we refer to the algorithm described in \cite{AlkoumiSchlenk2014}, where the higher dimensional generalization is straight-forward. For the case when the polytope is acute (what means that all dihedral angles are acute), we refer to the following result shown in \cite{AkopBal2015}: every length-minimizing closed billiard tarajectory on an acute convex polytope is maximally spanning, regular, and has exactly $n+1$ bouncing points. Therefore, for this case, it suffices to just consider the closed regular billiard trajectories with $n+1$ bouncing points.

The following is a sketch of the algorithm that produces length-minimizing closed regular billiard trajectories.

\begin{itemize}
\item[\small$\bullet$] For $m \in \{2,\ldots , n+1\}$ do:
\begin{itemize}
\item[\small$\bullet$] Choose $m$ pairwise different facets $F_1,\ldots ,F_m$ of $T$. For every $j \in \{1,\ldots ,m\}$ we let $u_j = n_T(p_{j+1})$ for some point $p_{j+1}$ in the relative interior of $F_{j+1}$. Note that $u_j$ does not depend on the choice of $p_{j+1}$.
\item[\small$\bullet$] Find $n_1,\ldots , n_m \in S^{n-1}\cap\spann\{u_1,...,u_m\}$ such that $n_{j+1} - n_{j} = -\mu_j u_{j}$ with $\mu_j > 0$ holds for every $j \in \{1,\ldots , m\}$.
\item[\small$\bullet$] Find $p_j \in F_j$ for $j\in \{1,\ldots,m\}$ such that $p_{j+1} - p_{j} = \lambda_j n_{j}$ with $\lambda_j > 0$ holds for every $j \in \{1,\ldots , m\}$.
\item[\small$\bullet$] Calculate the length of the closed polygonal curve and store it, if it is smaller than any such closed polygonal curve found so far.
\end{itemize}
\item[\small$\bullet$] Output the stored closed polygonal curve and its length.
\end{itemize} 
Proposition \ref{Prop:System} ensures that any closed polygonal curve found by the algo\-rithm is indeed a closed billiard trajectory. More precisely, if $p_1,\ldots ,p_m$ are the vertices of the closed polygonal curve, then \eqref{eq:System} is satisfied by construction. We will now examine this algorithm in more detail.

We start by letting $F_1,\ldots, F_m$ and $u_1,\ldots, u_m$ be as described above. Let $U$ be the $(n \times m)$-matrix that contains $u_1,\ldots, u_m$ as columns. Then it follows from Theorem \ref{Thm:RegularityResult1} that
\beq \textup{rk}(U)=m-1\label{eq:rankeucl}\eeq
is a necessary condition if we search for a length-minimizing closed billiard trajectory. Therefore, some choices of $F_1,\ldots, F_m$ can be discarded immediately. We note that \eqref{eq:rankeucl} also implies that
\begin{align*}
-\sum\limits_{j=1}^m \mu_j u_j = 0
\end{align*}
has up to scaling a unique solution $\mu_1,\ldots , \mu_m >0$. Consequently, there is (up to scaling) only one closed polygonal curve that can be constructed by using negative multiples of $u_1 \ldots , u_m$ in order.

To find suitable $n_1,\ldots ,n_m \in S^{n-1}\cap\spann\{u_1,...,u_m\}$, we let $\gamma = (\gamma_1,\ldots , \gamma_m)$ be this closed polygonal curve, encoded by its $m$ vertices. The task is now to scale (only using a positive scalar factor) and translate $\gamma$ such that the vertices of $\gamma$ lie on $S^{n-1}\cap\spann\{u_1,...,u_m\}$. We take $n_1,\ldots ,n_m$ as these vertices

Note that in the remainder of the algorithm it is required to form another closed polygonal curve using only positive multiples of $n_1,\ldots ,n_m$ in order, i.e.,
\begin{align}\label{eq:tcconfig}
\exists \lambda_1,\ldots , \lambda_m > 0 \colon \sum\limits_{j=1}^m \lambda_j n_j = 0.
\end{align}
If this property is true, then one says that $n_1,\ldots ,n_m$ are a \textit{totally cyclic vector configura\-tion}. Following \cite{Ziegler1995}, we can find an equivalent property by using Farkas' lemma. This property states that for every vector $v \in \mathbb{R}^n$, one of the following conditions hold:
\begin{itemize}
\item[(a)] $\inner{n_j}{v} < 0$, for some $j \in \{1,\ldots, m\}$,
\item[(b)] $\inner{n_j}{v} = 0$, for all $j \in \{1,\ldots, m\}$.
\end{itemize}
Hence, if $n_1,\ldots ,n_m$ are a totally cyclic vector configuration, it follows that
\begin{align}\label{unitPosition}
\forall\ v \in \mathbb{R}^n \ \exists\ j \in \{1,\ldots,m\}\ \colon \inner{v}{n_j} \leq 0.
\end{align}
This property is less restrictive than \eqref{eq:tcconfig}, but is still sufficient for the upcoming arguments.

While there might be multiple possibilities to scale and translate $\gamma$ such that its vertices lie on $S^{n-1}$, there is at most one possibility such that $\gamma$'s vertices $n_1,...,n_m$ on $S^{n-1}$ are a totally cyclic vector configuration (and this would also guarantee that the vertices of $\gamma$ lie not only on $S^{n-1}$, but on $S^{n-1}$ intersected with $\spann\{u_1,...,u_m\}$).

\bprop\label{Prop:biggestLine}
Let
\beqq \gamma = (\gamma_1,... ,\gamma_m) \in \mathbb{R}^n\times ... \times \mathbb{R}^n,\eeqq
\beqq S_\gamma = \left\{(\mu,t) \colon \mu\in \mathbb{R}_{\geq 0},\, t\in \mathbb{R}^n,\, ||\mu\gamma_j+t|| = 1,\, j\in \{ 1,...,m\}\right\} \eeqq
and
\beqq (\mu^\ast,t^\ast) \in S_\gamma\eeqq
be given. If $n_1,... ,n_m\in S^{n-1}$ satisfy \eqref{unitPosition}, where
\beqq n_j = \mu^\ast\gamma_j+t^\ast\eeqq
for every $j$, then
\beqq \mu^\ast = \max\{\mu \colon (\mu,t) \in S_\gamma \textup{ for some } t \in \mathbb{R}^n\}.\eeqq
\eprop

\bpf
We show that the existence of $(\widetilde{\mu},\widetilde{t}\,) \in S_\gamma$  with $\widetilde{\mu} > \mu^\ast$ yields a contradiction. For $j\in \{1,\ldots ,m\}$ we have
\begin{align*}
\widetilde{\mu} \gamma_j + \widetilde{t} = \widetilde{\mu} \frac{n_j - t^\ast }{\mu^\ast } + \widetilde{t} = \mu' n_j + t',
\end{align*}
where
\beqq \mu' = \widetilde{\mu}/\mu^\ast > 1 \; \text{ and } \; t' = \widetilde{t} - (\widetilde{\mu}/\mu^\ast) t^\ast.\eeqq
Because of $(\mu^\ast,t^\ast), (\widetilde{\mu},\widetilde{t}) \in S_\gamma$ we get
\begin{align*}
1  = ||\widetilde{\mu} \gamma_j + \widetilde{t}||^2 = ||\mu' n_j + t'||^2 &= {\mu'}^2 ||n_j||^2 +2\mu' \inner{n_j}{t'} + ||t'||^2\\ &= {\mu'}^2 + 2\mu' \inner{n_j}{t'} + ||t'||^2.
\end{align*}
Since $n_1,\ldots ,n_m$ satisfy \eqref{unitPosition}, there is some $i \in \{1,\ldots ,m\}$ such that $\inner{-t'}{n_i} \leq 0$. Note that this is obviously true if $t' = 0$. Therefore,
\begin{align*}
0 = {\mu'}^2 + 2\mu' \inner{n_i}{t'} + ||t'||^2 -1 \geq {\mu'}^2 + ||t'||^2 -1 > ||t'||^2 \geq 0,
\end{align*}
which is a contradiction.
\epf

We point out that the maximum in Proposition \ref{Prop:biggestLine} is indeed a maximum because the unit ball is compact and hence
\beqq \{\mu \colon (\mu,t) \in S_\gamma \textup{ for some } t \in \mathbb{R}^n\}\eeqq
is compact as well. Proposition \ref{Prop:biggestLine} implies that if the vertices of
\beqq \mu_1\gamma +t_1 \; \text{ and } \; \mu_2\gamma +t_2 \; \text{ with } \;(\mu_1,t_1), (\mu_2,t_2) \in S_\gamma\eeqq
are totally cyclic vector configurations (and hence satisfy \eqref{unitPosition}) then they are scaled by the same factor, i.e. $\mu_1 = \mu_2$. However, we also need $t_1 = t_2$ to make sure that there is only one suitable way to scale and translate $\gamma$. The next Proposition shows that this is indeed the case.

\bprop\label{Prop:uniqueDual}
Let
\beqq \gamma = (\gamma_1,\ldots ,\gamma_m) \in \mathbb{R}^n \times \ldots \times \mathbb{R}^n,\eeqq
$S_\gamma$ as in Proposition \ref{Prop:biggestLine}, and
\beqq (\mu_1,t_1),(\mu_2,t_2) \in S_\gamma.\eeqq
Further, let
\beqq n_1,\ldots ,n_m, n_1',\ldots n_m'\in S^{n-1}\eeqq
be defined by 
$$n_j = \mu_1\gamma_j + t_1 \, , \ n_j'=\mu_2\gamma_j + t_2$$
for every $j \in \{1,\ldots ,m\}$. If $n_1,\ldots,n_m$ satisfy \eqref{unitPosition}, then
\beqq (\mu_1,t_1)=(\mu_2,t_2),\eeqq
i.e., $n_j=n_j'$ for every $j \in \{1,\ldots ,m\}$.
\eprop

\bpf
From Proposition \ref{Prop:biggestLine} we have $\mu_1 = \mu_2$ and therefore for every $j\in \{1,\ldots, m\}$:
\begin{align*}
n_j' = \mu_1 \gamma_j + t_2 = \mu_1 \frac{n_j - t_1 }{\mu_1 } + t_2 = n_j - t_1 +t_2.
\end{align*}
Similar to the calculation in the proof of Proposition \ref{Prop:biggestLine}, we get
\begin{align*}
1 = || n_j' || = || n_j - t_1 + t_2 || = 1 + 2\inner{n_j}{t_2-t_1} +||t_2-t_1||^2.
\end{align*}
Hence, the term $\inner{n_j}{t_2-t_1}$ does not depend on $j$. Because of \eqref{unitPosition}, one can find $i,j \in \{1,\ldots , m\}$ such that 
$$\inner{n_i}{t_1-t_2} \leq 0 \; \text{ and }\; \inner{n_j}{t_1-t_2} \geq 0.$$ Since these terms don't depend on the index of $n$, we get
\beqq \inner{n_j}{t_1-t_2} = 0\eeqq
for every $j\in \{1,\ldots, m\}$. Further, we get
\begin{align*}
1 = 1 + ||t_2-t_1||^2,
\end{align*}
and so $t_1 = t_2$.
\epf

Assume we found some $n_1,\ldots, n_m\in S^{n-1}$ by scaling and translating the polygonal curve $\gamma$ as described above. Recall that, by construction of $\gamma$, we have
\beqq n_{j+1}-n_{j} = -\mu_j u_j\eeqq
for some $\mu_j > 0$ for every $j\in\{1,\ldots,m\}$. This implies
\begin{align*}
 1 = || n_{j+1} ||^2 = || n_{j} -\mu_j u_j  ||^2 &= || n_{j} || - 2 \mu_j \inner{n_{j}}{u_j} + \mu_j^2 ||u_j||^2 \\
 &= 1 - 2 \mu_j \inner{n_{j}}{u_j} + \mu_j^2
\end{align*}
and, consequently,
\beqq 0 = \mu_j (\mu_j - 2 \inner{n_{j}}{u_j}).\eeqq
Because $\mu_j$ is positive, we have
\beqq \mu_j = 2 \inner{n_{j}}{u_j}\eeqq
and
\begin{align*}
n_{j+1} = n_{j} -\mu_j u_j = n_{j} -2 \inner{n_{j}}{u_j} u_j
\end{align*}
for every $j \in \{1,\ldots ,m\}$. Therefore, if any of the vectors $n_1,\ldots, n_m$ is known, the remaining ones can be calculated easily. For what follows, we choose to search for $n_1$. This search can be carried out by a second-order cone program (SOCP).
An SOCP is a convex optimization problem, where one is looking for an element of the second-order cone
\beqq \mathcal{L}^{n+1} = \{(x,t) \in \mathbb{R}^{n}\times \mathbb{R} \colon ||x|| \leq t\}\eeqq
such that some linear constraints are satisfied and a linear objective function is optimized. An SOCP in standard form looks like this:
\begin{align*}
\textup{maximize }\hspace{2mm} & c^Tx + st\\
\textup{s.t. } & (x,t) \in \mathcal{L}^{n+1}\\
& a_i^Tx + r_it = b_i, \textup{ for } i \in I,
\end{align*}
where $c,a_i \in \mathbb{R}^n$, $s, r_i, b_i\in \mathbb{R}$ and $I$ is some index set. It is well known that SOCPs can be solved efficiently (see \cite{BenTalNemirovski2001}).

Before we state the SOCP, we need the following identity:
\begin{align*}
n_j = n_{j-1} -2 \inner{n_{j-1}}{u_{j-1}} u_{j-1} &= (I-2u_{j-1}u_{j-1}^T)n_{j-1} \\
&= (I-2u_{j-1}u_{j-1}^T)(I-2u_{j-2}u_{j-2}^T)n_{j-2} \\
&=\ldots = \left(\prod\limits_{i=1}^{j-1} I-2u_iu_i^T\right) n_1,
\end{align*}
where $I$ is the $n\times n$ identity matrix. This identity holds for $1\leq j \leq m$. Two types of constraints are necessary for the SOCP. First, we need to make sure that
\beqq n_1 = n_{m+1} \; \text{ with } \; n_{m+1}:= n_{m} -2 \inner{n_{m}}{u_m} u_m\eeqq
is satisfied:
\begin{align*}
 n_1 = n_{m} -2 \inner{n_{m}}{u_m} u_m = (I-2u_mu_m^T)n_{m} = \left(\prod\limits_{i=1}^m I-2u_iu_i^T\right) n_1.
\end{align*}
Second, we require
\begin{align*}
0 < \mu_j = 2 \inner{n_{j}}{u_j} = 2u_j^T\left(\prod\limits_{i=1}^{j-1} I-2u_iu_i^T\right) n_1
\end{align*}
for every $j \in \{1,\ldots,m\}$. Here, we can replace the strict inequality by $\leq$. The reason for this is the observation made earlier, that the polygonal curve $\gamma$ is unique up to scaling. So, if any solution of the SOCP yields $\mu_i =0$ for some $i$, then every solution yields $\mu_i =0$ and no closed regular billiard trajectory $p=(p_1,\ldots ,p_m)$ with $p_j\in F_j$ for $j \in \{1,\ldots ,m\}$ exists.

For the objective of the SOCP we note that 
\begin{align}\label{eq:pollinesize}
\sum\limits_{j=1}^m \mu_j = \sum\limits_{j=1}^m 2 \inner{n_{j}}{u_j} =  \sum\limits_{j=1}^m 2u_j^T \left( \prod\limits_{i=1}^{j-1} I-2u_iu_i^T \right)n_1.
\end{align}
The vectors $n_1,\ldots,n_m$ are required to be vertices of $\mu \gamma + t$ for some $(\mu,t)\in S_\gamma$. Proposition \ref{Prop:biggestLine} states that the only possible way for $n_1,\ldots,n_m$ to be a totally cyclic vector configuration is if $\mu$ is maximal. Hence, we would like to choose $n_1$ such that \eqref{eq:pollinesize} is as big as possible. Thus, we obtain the following SOCP:
\begin{align*}
\textup{maximize }\hspace{2mm} & \sum\limits_{j=1}^m 2u_j^T \left( \prod\limits_{i=1}^{j-1} I-2u_iu_i^T \right)x\\
\textup{s.t. }\hspace{2mm} & x\in \mathbb{R}^n,\ || x || \leq 1 \\
& 2u_j^T\left(\prod\limits_{i=1}^{j-1} I-2u_iu_i^T\right) x \geq 0, \ \forall \ j\in \{1,\ldots ,m\}\\
& \left(\left(\prod\limits_{i=1}^m I-2u_iu_i^T\right)-I\right) x = 0.
\end{align*}
If this SOCP has an optimal solution $x^\ast$, we pick $n_1=x^\ast$. It is easy to see that this way $n_1$ lies on $S^{n-1}$. If we assume
\beqq ||x^\ast|| = ||n_1|| <1,\eeqq
then $x^\ast/||x^\ast||$ would be feasible as well, but would have a greater objective value (unless the optimal value is $0$ in which case $\mu_j=0$ for every $j$. But then, there is no closed regular billiard trajectory for the given choice of facets). This contradicts the optimality of $x^\ast$.

We now have a way to find $n_1$ and, hence, also $n_2,\ldots ,n_m$. These vectors are unique by Proposition \ref{Prop:uniqueDual}. We proceed in a similar fashion as before in order to find the bouncing points $p_1,\ldots ,p_m$ of a closed regular billiard trajectory. \eqref{eq:System} states that
\beq p_{j+1} - p_j = \lambda_j n_j\label{eq:normalvectorspanningeuclr}\eeq
needs to hold for every $j \in \{1,\ldots , m\}$, where $\lambda_j > 0$. Thus, just like $u_1,\ldots ,u_m$, we have that $n_1,\ldots ,n_m$ are not linearly independent:
\begin{align*}
0 = \sum\limits_{j=1}^m (p_{j+1} - p_j) = \sum\limits_{j=1}^m \lambda_j n_j.
\end{align*}
We would like to have that the $(n\times m)$-matrix $(n_1,\ldots , n_m)$ has rank $m-1$. Then, the closed polygonal curve which can be constructed by positive multiples of $n_1,\ldots , n_m$ is unique up to scaling. Indeed, considering \eqref{eq:normalvectorspanningeuclr}, this is the case as direct consequence of Theorem \ref{Thm:RegularityResult1} and an argument from within the proof of Proposition \ref{Prop:normalvectorsspanning}.

For what follows, we let $\xi = (\xi_1,\ldots ,\xi_m)$ be such a closed polygonal curve, i.e.,
\beqq \xi_{j+1} -\xi_{j} = \lambda'_j n_j\eeqq
with $\lambda'_j>0$ for $j\in \{1,\ldots ,m\}$. The task is now to find $\lambda > 0$ and $s \in \mathbb{R}^n$ such that
\beqq \lambda \xi_j + s \in F_j\eeqq
for every $j$. This can be done via a linear program (LP).

Just like an SOCP, an LP is a convex optimization with a linear objective function and linear constraints. As opposed to an SOCP, where the goal is to find an optimal member of $\mathcal{L}^{n+1}$, an LP asks for an optimal member of the non-negative orthant $\mathbb{R}^n_{\geq 0}$. An LP in standard form looks as follows:
\begin{align*}
\textup{maximize }\hspace{2mm} & c^Tx\\
\textup{s.t. } & x \in \mathbb{R}^n_{\geq 0}\\
& Ax = b,
\end{align*}
where $c \in \mathbb{R}^n, b \in \mathbb{R}^r, A \in \mathbb{R}^{r\times n}$. Here, $r$ is the number of linear constraints. Linear programming is a special case of second-order cone programming. In particular, it is not surprising that LPs can be solved efficiently as well (see \cite{Schrijver1998}).

For $j \in \{1,\ldots ,m\}$ let $$H_j = \{x \in \mathbb{R}^n \colon \inner{u_j}{x} = b_j \}$$ be the unique supporting hyperplane of $T$ which contains $F_j$. Then we require
\begin{align*}
\lambda \xi_j + s \in F_j & \Longleftrightarrow\hspace*{1mm} \lambda \xi_j + s \in T \cap H_j\\
& \Longleftrightarrow\hspace*{1mm} \inner{\lambda \xi_j + s}{u_j} = b_j,\ \inner{\lambda \xi_j + s}{u_i} \leq b_i,\ \forall \ i\neq j\\
& \Longleftrightarrow\hspace*{1mm} \begin{pmatrix}
\xi_j^Tu_j & u_j^T 
\end{pmatrix} \begin{pmatrix}
\lambda\\
s
\end{pmatrix} = b_j,\ \begin{pmatrix}
\xi_j^Tu_i & u_i^T 
\end{pmatrix} \begin{pmatrix}
\lambda\\
s
\end{pmatrix} \leq b_i,\ \forall \ i\neq j.
\end{align*}
Additionally, we would like to make sure that (if possible) the resulting closed billiard trajectory
\beqq p = (p_1,\ldots ,p_m) \; \text{ with } \; p_j = \lambda \xi_j + s\eeqq
is regular. Thus, we maximize the inclusion minimal slack of the form
\beqq b_i - (\lambda \xi^T_ju_i + u_i^Ts) \; \text{ with } \; 1\leq i,j\leq m \; \text{ and } \; i\neq j.\eeqq

This leads us to the following LP with $n+2$ variables:

\begin{align*}
\textup{maximize }\hspace{2mm} & \rho\\
\textup{s.t. }\hspace{2mm} & \rho,\lambda\in \mathbb{R}_{\geq 0},\ s \in \mathbb{R}^n \\
& \begin{pmatrix}
\xi_j^Tu_j & \,u_j^T 
\end{pmatrix} \begin{pmatrix}
\lambda\\
s
\end{pmatrix} = b_j, \ \forall \ j \in \{1,\ldots ,m\}\\
& \rho \leq b_i-\begin{pmatrix}
\xi_j^Tu_i & \,u_i^T 
\end{pmatrix} \begin{pmatrix}
\lambda\\
s
\end{pmatrix},\ \forall i,j\in \{1,\ldots ,m\}, \ i\neq j.
\end{align*}
Note that any solution of the LP satisfies $\lambda > 0$. If otherwise $\lambda = 0$, then the equality constraints state that $s\in H_j$. This means that all supporting hyperplanes which contain a face of $T$ intersect at a common point. This is not possible. If a solution exists, we let
\beqq p_j = \lambda\xi_j + s\eeqq
for $j\in \{1,\ldots ,m\}$. If the optimal solution satisfies $\rho > 0$, then $p = (p_1,\ldots ,p_m)$ is a closed regular billiard trajectory which is potentially length-minimizing (meaning that it satisfies the condition in Theorem \ref{Thm:RegularityResult1}). To ensure that we don't miss any potentially length-minimizing closed regular billiard trajectories, we are going to prove that it suffices to find one such closed billiard trajectory per choice of $F_1,\ldots, F_m$.

\bprop\label{Prop:SameLength}
Let $F_1,\ldots ,F_m$ be facets of some full-dimensional polytope $T$ and let $p_j,p_j' \in F_j$ for every $j$. Assume
\beqq p = (p_1,\ldots ,p_m)\; \text{ and } \; p' = (p'_1,\ldots ,p'_m)\eeqq
are closed regular billiard trajectories in $T$. Then,
\beqq \ell(p) = \ell(p').\eeqq
\eprop

\bpf
According to Proposition \ref{Prop:System}, there are $n_1,\ldots ,n_m, n_1',\ldots ,n_m'$ such that:
\begin{align*}
\begin{cases}p_{j+1}-p_j = \lambda_j n_j, \; \lambda_j>0, \\ n_{j+1}-n_j \in -N_T(p_{j+1}), \end{cases} \textup{ and }\;\; \begin{cases}p'_{j+1}-p'_j = \lambda'_j n'_j, \; \lambda'_j>0, \\ n'_{j+1}-n'_j \in -N_T(p'_{j+1}). \end{cases}
\end{align*}
As has been discussed earlier in this section, $n_1,\ldots ,n_m$ are unique and only depend on $n_T(p_1),\ldots , n_T(p_m)$. Therefore, we have $n_j' = n_j$ for $j \in \{1,\ldots , m\}$. Now, one has natural projection maps from $H_1$ to $H_2$ along $n_1$, from $H_2$ to $H_3$ along $n_2$, and so on, where we denote by $H_j$ the $T$-supporting hyperplane that contains $F_j$. Their composition $H_1\rightarrow H_1$ is an affine map that decreases the lengths of vectors in $H_1\cap\spann\{n_i\}$; hence, it cannot have more than one fixed point.
%Furthermore, we can use the fact that $||n_j|| = 1$ to get:
%\begin{align*}
%\ell(p) &= \sum\limits_{j=1}^m ||p_{j+1} -p_j|| = \sum\limits_{j=1}^m ||\lambda_j n_j|| = \sum\limits_{j=1}^m \lambda_j ||n_j|| = \sum\limits_{j=1}^m \lambda_j \inner{n_j}{n_j} \\&= \sum\limits_{j=1}^m \inner{\lambda_j n_j}{n_j} = \sum\limits_{j=1}^m \inner{p_{j+1} -p_j}{n_j}.
%\end{align*}
%Similarly, we get: $$\ell(p') = \sum_j \inner{p'_{j+1} -p'_j}{n_j}.$$ This implies:
%\begin{align*}
%\ell(p)-\ell(p') &= \sum\limits_{j=1}^m \inner{p_{j+1} -p_j}{n_j} - \inner{p'_{j+1} -p'_j}{n_j} \\&= \sum\limits_{j=1}^m \inner{p_{j+1} -p_{j+1}'}{n_j} + \inner{p'_{j} -p_j}{n_j}\\
%&= \sum\limits_{j=1}^m \inner{p_{j+1} -p_{j+1}'}{n_j} + \sum\limits_{j=1}^m\inner{p'_{j+1} -p_{j+1}}{n_{j+1}}\\
%&= \sum\limits_{j=1}^m \inner{p'_{j+1} -p_{j+1}}{n_{j+1}-n_j}.
%\end{align*}
%As before, we let $$H_{j+1} = \{x\in \mathbb{R}^n \colon \inner{n_T(p_{j+1})}{x} = b_{j+1}\}$$ be the supporting hyperplane of $T$ which contains $F_{j+1}$. On the one hand we have
%\begin{align*}
%\inner{n_T(p_{j+1})}{p'_{j+1} -p_{j+1}} &= \inner{n_T(p_{j+1})}{p'_{j+1}} - \inner{n_T(p_{j+1})}{p_{j+1}}\\
%& =  b_{j+1} - b_{j+1} = 0.
%\end{align*}
%On the other hand $n_{j+1}-n_j \in -N_T(p_{j+1})$. Thus $n_{j+1}-n_j$ is a multiple of $n_T(p_{j+1})$. This yields $\ell(p)-\ell(p') = 0$.
\epf

For each choice of $F_1,\ldots ,F_m$, the algorithm needs to solve the following tasks: Calculate the rank of an $(n\times m)$-matrix, solve an SOCP with $n+m+1$ constraints and $n+1$ variables, solve an $(n+1)\times m$ system of linear equations, and solve an LP with $mf$ constraints and $n+2$ variables, where $f$ is the number of facets of $T$. All these tasks are solvable in polynomial time (with respect to the dimension $n$ and the number of facets of $T$); see for instance \cite{Nesterov1994}. However, there are
\beqq \sum_{j=2}^{n+1} \binom{f}{j}j!\eeqq
possibilities to choose at least $2$ but at most $n+1$ facets, respecting their order. We can slightly improve this number since a cyclic shift of the chosen facets $F_1,\ldots ,F_m$ will yield a similar (but shifted) result. This leaves us with
\beqq \sum_{j=2}^{n+1} \binom{f}{j}(j-1)!\eeqq
possibilities. The calculations for each of these possibilities are independent of each other. Therefore, we utilize parallel computing to accelerate the algorithm.

In Table \ref{table:times}, the running time of the algorithm can be seen. Each time, the billiard table is a polytope $T$ with dimension $2$, $3$, or $4$ which we generated in the following way: First, we chose some normally distributed random vectors. We scaled each of these vectors by some scalar between $1$ and $3$ (we decreased the length of this range if the amount of random vectors became too large). Afterwards, we received $T$ as the convex hull of these vectors. Instead of a total running time, the table shows the time needed to compute a length-minimizing closed regular billiard trajectory with $2$, $3$, $4$, and $5$ bouncing points, respectively. The table suggests that the calculations for $m$ bouncing points with $m < d+1$ terminate very quickly. The reason for this is that many iterations are cancelled early when the rank of $(u_1,\ldots , u_m)$ is checked. All calculations have been done on a Dell Latitude E6530 laptop with Intel Core i7-3520M processor, 2.9 GHz (capable of running four threads). The algorithm has been implemented in Python and mainly utilizes the NumPy library. The LPs and SOCPs are solved via the software Mosek. The algorithm and a detailed description on how to choose the input is available on the website \href{www.github.com/S-Krupp/EHZ-capacity-of-polytopes}{www.github.com/S-Krupp/EHZ-capacity-of-polytopes} (cf.\;also the description in \cite{Krupp2021}).

\begin{table}\label{Table:Times}
\centering
\caption{Running times of the algorithm outlined above. The billiard table $T$ is a polytope. The first two columns contain the number of facets and dimension of $T$. The last four columns contain the running time for 2, 3, 4 and 5 bouncing points in seconds.}\label{table:times}
\begin{tabular}[h]{c|c|c|c|c|c}
\# facets & dim & time 2 bp. & time 3 bp. & time 4 bp. & time 5 bp.\\
\hline
 10 & 2 & 0.00831 &    0.30746 & - & - \\ 
 15 & 2 & 0.00973 &    0.77826 & - & - \\ 
 20 & 2 & 0.01228 &    1.90130 & - & - \\ 
 25 & 2 & 0.01535 &    3.68928 & - & - \\ 
 30 & 2 & 0.01885 &    6.87248 & - & - \\ 
 35 & 2 & 0.02242 &   10.98917 & - & - \\ 
 40 & 2 & 0.02611 &   17.82252 & - & - \\ 
 45 & 2 & 0.02676 &   26.24746 & - & - \\ 
 50 & 2 & 0.03756 &   37.58889 & - & - \\ 
 60 & 2 & 0.04348 &   68.84613 & - & - \\ 
 70 & 2 & 0.07126 &  121.99939 & - & - \\ 
 80 & 2 & 0.08522 &  181.44979 & - & - \\ 
 90 & 2 & 0.09377 &  276.15113 & - & - \\ 
100 & 2 & 0.12422 &  390.50119 & - & - \\ 
110 & 2 & 0.14459 &  516.93047 & - & - \\ 
120 & 2 & 0.16201 &  706.91467 & - & - \\ 
130 & 2 & 0.16660 &  887.57707 & - & - \\ 
140 & 2 & 0.23777 & 1145.28408 & - & - \\ 
150 & 2 & 0.26617 & 1400.09367 & - & - \\ 
\hline 
 14 & 3 & 0.00983 & 0.02658 &    2.78469 & - \\ 
 20 & 3 & 0.01315 & 0.05639 &   12.75729 & - \\ 
 24 & 3 & 0.01319 & 0.10031 &   25.36455 & - \\ 
 30 & 3 & 0.01869 & 0.16165 &   69.23203 & - \\ 
 34 & 3 & 0.01990 & 0.30546 &  121.10618 & - \\ 
 40 & 3 & 0.02877 & 0.43128 &  281.39158 & - \\ 
 44 & 3 & 0.02691 & 0.50671 &  456.25295 & - \\ 
 50 & 3 & 0.03799 & 0.76186 &  755.02158 & - \\ 
 54 & 3 & 0.04065 & 0.96421 & 1091.77615 & - \\ 
 60 & 3 & 0.04458 & 1.32361 & 1646.65092 & - \\ 
 64 & 3 & 0.04991 & 1.61337 & 2158.03637 & - \\ 
 70 & 3 & 0.07663 & 2.09306 & 2849.52804 & - \\ 
\hline
 11 & 4 & 0.00991 & 0.01654 &  0.05954 &    2.36219 \\ 
 15 & 4 & 0.00870 & 0.02263 &  0.17406 &   19.43698 \\ 
 20 & 4 & 0.01176 & 0.04806 &  0.56930 &   54.60817 \\ 
 25 & 4 & 0.01730 & 0.10708 &  1.58835 &  245.41436 \\ 
 30 & 4 & 0.02492 & 0.24002 &  4.56021 &  961.83634 \\ 
 35 & 4 & 0.02928 & 0.36108 &  9.92964 & 2171.25146 \\ 
 40 & 4 & 0.03232 & 0.50996 & 19.76087 & 4201.35654 
\end{tabular}
\end{table}

\subsection{Constructing closed non-regular billiard trajectories}

We already noted that it is much more complicated to construct closed non-regular billiard trajectories (with more than two bouncing points and on higher dimensional tables). From Theorem \ref{Thm:RegularityResult1} we know that the length-minimizing closed billiard trajectories are maximally spanning and have at most $n+1$ bouncing points, but they do not have to be relatively regular, i.e., they do not have to be regular within the inclusion minimal affine section that contains themselves. This makes it difficult to reduce the problem of searching for length-minimizing closed billiard trajectories to the problem of searching for suitable closed billiard trajectories within affine sections of the billiard table.

One way to overcome this problem is to set a condition on the billiard tables that are permitted such that the possibility of having length-minimizing closed billiard trajectories which are not relatively regular is excluded. With a look at the proof of Theorem \ref{Thm:RegularityResult1}, this can be achieved by requiring the convex bodies $T\subset\R^n$ to satisfy the following condition: for every affine section $T\cap V$ of $T$ that contains a set of points on $\partial T$ that is in $F(T)$, we have
\beq T\subset (T\cap V)\oplus W,\label{eq:conditionrelativelyregular}\eeq
where $W$ is the $(n-\dim V)$-dimensional linear subspace orthogonal to $V$\footnote{We remark that one could weaken the condition by requiring that for every affine section $T\cap V$ of $T$ that contains a (length-minimizing) closed billiard trajectory we have \eqref{eq:conditionrelativelyregular}. However, this condition is not so easy to check apriori, because we are just searching for (length-minimizing) closed billiard trajectories.}. We call the set of convex bodies in $\R^n$ satisfying \eqref{eq:conditionrelativelyregular} $\mathcal{C}^{c}(\R^n)$ and conclude the following corollary:

\bcor\label{Thm:relativelyregular}
Let $T\in \mathcal{C}^c(\R^n)$ be a billiard table and $p=(p_1,...,p_m)$ a length-minimizing closed billiard trajectory on $T$. Further, let $U\subseteq\R^n$ be the convex cone spanned by the normal vectors related to the billiard reflection rule and let $V\subseteq\R^n$ be the affine subspace such that $T\cap V$ is the inclusion minimal affine section of $T$ that contains $p$. Then, one has $U=V_0$, where $V_0$ is the linear subspace underlying $V$, and $p$ is maximally spanning and relatively regular.
\ecor

One can easily see from the last part of the proof of Theorem \ref{Thm:RegularityResult1} (when arguing that the constructed closed polygonal curve $\widetilde{p}$ is in $F(T)$) that the condition $T\in\mathcal{C}^c(\R^n)$ allows one to replace \eqref{eq:RegularityResult1}, i.e., 
\beqq \dim(N_T(p_j)\cap V_0)=1 \quad \forall j\in\{1,...,m\}\eeqq
by
\beqq \dim N_{T\cap V}(p_j)=1\quad \forall j\in\{1,...,m\}.\eeqq

Let $P\in\mathcal{C}^c(\R^n)$ now be a convex polytope. Now, using Corollary \ref{Thm:relativelyregular}, our aim is to construct closed billiard trajectories on $P$ which come into question to minimize the length. More precisely, we construct maximally spanning and relatively regular closed billiard trajectories on $P$ with $j$ bouncing points for $j\in\{3,...,n+1\}$. 

For that, we first introduce the following useful characterization of affine sections of $P$: We call an affine section $P\cap V$ of $P$ \textit{essential} if
\beqq P\cap V \in F(P).\eeqq
By Lemma \ref{Lem:Bezdek}(ii), this can be equivalently expressed by: $P\cap V$ is essential if there are $P$-supporting half-spaces of $\R^n$ through points on $\partial (P\cap V)$ such that its intersection is nearly bounded.

Expressed in pseudo code, we then do the following steps for every essential $(j-1)$-dimensional affine section $P\cap V$ of $P$ (an explanation follows--also concerning whether this produces an endless loop):
\begin{itemize}
\item[(a)] Choose $j$ facets $F_1,...,F_j$ of $P\cap V$ (considering the order) such that the cone spanned by the associated outer normal unit vectors $n_{F_1},...,n_{F_j}$ (within $V$) is $V_0$.
\item[(b)] Construct the uniquely determined (up to scaling and translation) closed polygonal curve $(\gamma_1,...,\gamma_j)$, where the $\gamma_{i+1}-\gamma_i$ are given by negative multiples of the $n_{F_i}$ for all $i\in \{1,...,j\}$.
\item[(c)] If possible: Find $\lambda >0$ and $c\in \R^n$ such that
\beqq \lambda\{\gamma_1,...,\gamma_j\}+c\subset S^{n-1}\cap V_0.\eeqq
Otherwise: If possible: Go back to step (a) and start with a choice not yet made. Otherwise: End. 
\item[(d)] Define unit vectors $n_1,...,n_j$ by
\beqq n_i:=\lambda\gamma_i+c \quad \forall i\in\{1,...,j\}.\eeqq
Construct the uniquely determined (up to scaling and translation) closed polygonal curve $(\xi_1,...,\xi_j)$, where the $\xi_{i+1}-\xi_i$ are given by positive multiples of the $n_{i+1}$ for all $i\in \{1,...,j\}$. Otherwise: If possible: Go back to step (a) and start with a choice not yet made. Otherwise: End.
\item[(e)] If possible: Find $\mu >0$ and $e\in\R^n$ such that
\beqq \mu\xi_i+e\in \mathring{F}_i\quad \forall i\in\{1,...,j\}.\eeqq
Otherwise: If possible: Go back to step (a) and start with a choice not yet made. Otherwise: End.
\item[(f)] We define a closed polygonal curve $p=(p_1,...,p_j)$ by
\beqq p_i:=\mu\xi_i+e\quad \forall i\in\{1,...,j\}.\eeqq
By construction: $p$ is a maximally spanning and relatively regular closed billiard trajectory on $P$ with $j$ bouncing points. Add $p$ to the set $B_{j}(P)$.
\item[(g)] If possible: Go back to step (a) and start with a choice not yet made. Otherwise: End.
\end{itemize}
For every $j\in \{2,...,n+1\}$ the set $B_{j}(P)$ contains all maximally spanning and relatively regular closed billiard trajectories on $P$ with $j$ bouncing points whose length can be easily calculated by \eqref{eq:length}. The final task is to find at least one minimizer in every $B_j(P)$, and consequently one global minimizer.

Let us sketchily explain the idea behind the above steps: The main idea is to use Proposition \ref{Prop:System}. More precisely, for every $j\in\{3,...,n+1\}$, we first select $j$ facets of each essential $(j-1)$-dimensional affine section of $P$ which by Proposition \ref{Prop:normalvectorsspanning} (which is also reflected within Theorem \ref{Thm:RegularityResult1} and Corollary \ref{Thm:relativelyregular}) have to satisfy $U=V_0$, where $U$ is the convex cone spanned by the corresponding normal vectors, and check whether the choice of the facets (considering their order) makes the existence of closed billiard trajectories with bouncing points in the interiors of these facets possible. We check this implicitly by first using the second condition in \eqref{eq:System} to determine the unit vectors in \eqref{eq:System}. By a geometrical argument, we then draw conclusions from the first condition in \eqref{eq:System} as to whether the above mentioned closed billiard trajectories can exist. All of this can be done for each choice of $j$ facets (considering their order) of each essential $(j-1)$-dimensional affine section of $P$ so that after a finite number of steps we have constructed all regular closed billiard trajectories with $j$ bouncing points within each essential $(j-1)$-dimensional affine section of $P$.

Starting with the selection of essential $(j-1)$-dimensional affine sections of $P$, stems from Lemma \ref{Lem:Bezdek}(iii), i.e., from the fact that every closed billiard trajectory on $P$ is in $F(P)$. We, however, note that there can be an infinite number of essential affine sections of $P$. Nevertheless, since $P$ is a convex polytope, from the qualitative perspective, there are only finitely many different essential affine sections. With a clever implementation, one should be able to avoid qualitatively-equal choices of essential affine sections. This would prevent the algorithm to produce endless loops.

Let us now turn to the explanation of the individual steps:

Ad (a): If there is a maximally spanning and relatively regular closed billiard trajectory with $j$ bouncing points contained in an essential affine section $P\cap V$, then we know together with Proposition \ref{Prop:normalvectorsspanning}, i.e., with $U=V_0$, that the bouncing points lie in the interiors of $j$ different facets of $P\cap V$ while the cone spanned by the corresponding outer unit normal vectors is $V_0$. This implies by Proposition \ref{Prop:noteinproof} that the intersection of the half-spaces which contain $P$ and support these facets is nearly bounded in $\R^n$. By Lemma \ref{Lem:Bezdek}(ii), this ensures that a closed polygonal curve with vertices on these facets cannot be translated into the interior of $P$ (which by \ref{Lem:Bezdek}(iii) is a necessary condition for closed billiard trajectories on $P$). This justifies to start with essential affine sections of $P$ as necessary condition for the existence of closed billiard trajectories which come into question to minimize the length.

Ad (b): Since the cone spanned by $n_{F_1},...,n_{F_j}$ is $(j-1)$-dimensional, solving a system of linear equations yields a uniquely determined (up to scaling and translation) $j$-tuple $(\gamma_1,...,\gamma_j)\in(\R^n)^j$ and a uniquely determined (up to scaling) $j$-tuple $(\alpha_1,...,\alpha_j)\in (\R_{<0})^j$ satisfying
\beqq \gamma_{i+1}-\gamma_i=\alpha_i n_{F_i}\eeqq
for all $i\in\{1,...,j\}$. We understand the $j$-tuple $(\gamma_1,...,\gamma_j)$ as a closed polygonal curve.

Ad (c): There is at most one combination $(\lambda,c)\in \R_{>0}\times \R^n$ such that
\beqq \lambda\{\gamma_1,...,\gamma_j\}+c\subset S^{n-1}\cap V_0.\eeqq
This is guaranteed by the fact that the underlying linear subspace $V_0$ of $V$ is unique.

Ad (b) \& (c): If we define unit vectors $n_1,...,n_j$ by
\beqq n_i:=\lambda\gamma_i+c\eeqq
for all $i\in\{1,...,j\}$ while it is the case that the condition in step (c) holds for $(\lambda,c)\in\R_{>0}\times \R^n$, then these unit vectors are candidates to be the unit vectors in \eqref{eq:System} associated to a possibly existing closed billiard trajectory on $P$ with bouncing points in the interiors of $F_1,...,F_j$. However, the second condition in \eqref{eq:System} is already guaranteed. In what follows, it has to be checked whether there really is an associated closed billiard trajectory on $P$ with bouncing points in the interiors of $F_1,...,F_j$.

Ad (d): The cone spanned by the unit vectors $n_1,...,n_j$ is $V_0$. Therefore, as in step (b), solving a system of linear equations yields a uniquely determined (up to scaling und translation) $j$-tuple $(\xi_1,...,\xi_j)\in(\R^n)^j$ and a uniquely determined (up to scaling) $j$-tuple $(\beta_1,...,\beta_j)\in (\R_{>0})^j$ satisfying
\beqq \xi_{i+1}-\xi=\beta_{i+1} n_{i+1}\eeqq
for all $i\in\{1,...,j\}$. We understand the $j$-tuple $(\xi_1,...,\xi_j)$ as a closed polygonal curve.

Ad (e) \& (f): Similar to step (c), there is at most one combination $(\mu,e)\in \R_{>0}\times \R^n$ such that
\beqq \mu\{\xi_1,...,\xi_j\}+e\subset \partial (P\cap V).\eeqq
By checking whether
\beqq \mu \xi_i +e \in \mathring{F}_i \quad \forall i\in\{1,...,j\},\eeqq
we make sure that the closed polygonal curve $p=(p_1,...,p_j)$ defined by
\beqq p_i:=\mu \xi_i +e\eeqq
for all $i\in\{1,...,j\}$ has its vertices in the interiors of the facets $F_1,...,F_j$. The closed polygonal curve is a closed billiard trajectory on $P$ since $p=(p_1,...,p_j)$ is fulfilling \eqref{eq:System} together with the $j$-tuple of unit vectors $(n_2,...,n_j,n_1)$ and
\beqq \lambda_i:=\mu\beta_i>0,\; \mu_i:=\lambda \alpha_i<0 \quad \forall i\in\{1,...,j\}.\eeqq
Indeed, we have
\beqq \begin{cases} p_{i+1}-p_i=(\mu \xi_{i+1}+e)-(\mu\xi_i +e)=\mu(\xi_{i+1}-\xi_i)=\mu\beta_{i+1} n_{i+1}=\lambda_{i+1} n_{i+1},\\
n_{i+1}-n_i=(\lambda\gamma_{i+1}+c)-(\lambda\gamma_i+c)=\lambda (\gamma_{i+1}-\gamma_i)=\lambda\alpha_i n_{F_i}\in -N_P(p_i).\end{cases}
\eeqq
Additionally, $p$ is maximally spanning and relatively regular by construction. Therefore, $p$ is a maximally spanning and relatively regular closed billiard trajectory on $P$ which is contained in the essential affine section $P\cap V$.

Tasks we left for further research are the theoretical clarification of the individual steps as in Section \ref{Subsec:constructingregular}, the implementation of the algorithm (which includes to develop a technique in order to find only finitely many qualitatively-different essential affine sections), and a criteria for deciding whether a convex body/polytope is in $\mathcal{C}^c(\R^n)$. This last task seems to be challenging since, obviously, not every convex body can be approximated by convex polytopes in $\mathcal{C}^c(\R^n)$. For instance, one easily checks that for Example (C) in Section \ref{Sec:Examples} this is not possible. It remains the question how big $\mathcal{C}^c(\R^n)$ is in $\mathcal{C}(\R^n)$.

\section*{Acknowledgement}
This research was carried out under the supervision of Alberto Abbondandolo and Frank Vallentin. The authors are thankful to the supervisors' support and thank  Alexey Balitskiy for his remarks on a previous draft of this paper. The authors also thank the unknown referee for the useful comments.

%\bibliographystyle{plain}
%\bibliographystyle{alpha}
%\bibliography{Billiards}

\medskip

\section*{Daniel Rudolf, Ruhr-Universit\"at Bochum, Fakult\"at f\"ur Mathematik, Universit\"atsstrasse 150, D-44801 Bochum, Germany.}
\center{E-mail address: daniel.rudolf@ruhr-uni-bochum.de}

\section*{Stefan Krupp, Universit\"at zu K\"oln, Mathematisches Institut, Weyertal 86-90, D-50931 K\"oln, Germany.}
\center{E-mail address: krupp@math.uni-koeln.de}

\end{document}